\documentclass[11pt]{article}
\usepackage{latexsym,amssymb,amsmath}
\def\dstar{\delta_{*}}

\def\ff#1{\mathbb F\mathbb F^{#1}}
 \def\fg{\mathfrak g}
\def\fh{\mathfrak h} \def\fp{\mathfrak p}

\def\inv{{}^{-1}}

\def\La#1{\Lambda^{#1}}

\def\ot{\!\otimes\!}

\def\pp#1{\mathbb P^{#1}}
\def\ppp{\mathbb P}

\def\ra{\rightarrow}

\def\tbase{{\rm Base}\,}

\def\tdim{{\rm dim}\,}

\def\trank{{\rm rank}\,}

\def\ttrace{{\rm tr}}

\def\th{\theta}

\def\up#1{{}^{({#1})}}
\def\upperp{{}^{\perp}}

\def\ww{\wedge}

\def\a{{\alpha}}
\def\b{{\beta}}
\def\c{{\gamma}}
\def\d{{\delta}}

\def\l{{\lambda}}

\def\o{{\omega}}

\def\s{{\sigma}}
\def\t{{\tau}}
\def\th{{\theta}}

\def\ZZ{\mathbb Z}
\def\CC{\mathbb C}

\def\PP{\mathbb P}
\def\FF{\mathbb F}
\def\SS{\mathbb S}
\def\RR{\mathbb R}
\def\QQ{\mathbb Q}

\def\11{{\mathord{\bf 1\/}}}
\def\ot{{\mathord{\otimes }}}

 \def\we{{\mathord{\Wedge }}}

\def\ra{{\mathord{\;\rightarrow\;}}}

\def\com{{\rm com}}
\def\trace{{\rm trace}}

\def\cZ{{\cal Z}}
\def\CC{\mathbb C}
\def\RR{\mathbb R}
\def\HH{\mathbb H}
\def\AA{{\mathbb A}}
\def\BB{{\mathbb B}}
\def\OO{\mathbb O}
\def\BC{\mathbb C}
\def\BR{\mathbb R}
\def\BH{\mathbb H}
\def\BA{{\mathbb A}}
\def\BB{{\mathbb B}}
\def\BO{\mathbb O}
\def\BP{\mathbb P}
\def\BF{\mathbb F}
\def\BS{\mathbb S}

\def\ZZ{\mathbb Z}
\def\SS{\mathbb S}

\def\11{\mathbf 1}
\def\PP{\mathbb P}
\def\QQ{\mathbb Q}
\def\FF{\mathbb F}
\def\JA{{\cal J}_3(\AA)}
\def\JB{{\cal J}_3(\BB)}

\def\fh{{\mathfrak h}}
\def\fsl{{\mathfrak {sl}}}
\def\fsp{{\mathfrak {sp}}}

\def\fso{{\mathfrak {so}}}
\def\fs{{\mathfrak {s}}}
\def\fe{{\mathfrak e}}
\def\ff{{\mathfrak f}}

\def\fg{{\mathfrak g}}

\def\fp{{\mathfrak p}}

\def\fl{{\mathfrak l}}
\def\l{\lambda}
\def\a{\alpha}
\def\o{\omega}
\def\om{\omega}

\def\O{\Omega}
\def\b{\beta}

\def\s{\sigma}

\def\d{\delta}

\def\ot{{\mathord{\otimes }\,}}
\def\op{{\mathord{\oplus }\,}}

\def\ra{{\mathord{\;\rightarrow\;}}}
\def\we{{\mathord{{\scriptstyle \wedge}}}} \def\tr{{\rm trace}\;}
\newcommand\qed{{\hspace*{\fill} $\Box$}} 
\newcommand\proof{{\noindent {\em Proof}.}\hspace{2mm}}
\newcommand\proofs{{\noindent {\em Proofs}.}\hspace{2mm}}
\newcommand\rem{{\noindent {\em Remark}.}\hspace{2mm}}

\newtheorem{theo}{Theorem}[section]
\newtheorem{coro}[theo]{Corollary}
\newtheorem{lemm}[theo]{Lemma}
\newtheorem{prop}[theo]{Proposition}

\begin{document}

\title{The projective geometry of Freudenthal's magic square} \author{J.M.
Landsberg and L. Manivel}
\date{}
\maketitle

\begin{abstract} We connect the algebraic geometry and representation
theory associated to Freudenthal's magic square. We give 
unified geometric descriptions of several classes of orbit closures,
describing their hyperplane sections and desingularizations, and
interpreting them in terms of composition algebras. In particular, we
show how a class of invariant quartic polynomials can be viewed 
as generalizations of
the classical discriminant of a cubic polynomial.\end{abstract} 

{\small \tableofcontents}

\section{Introduction}
\subsection{The magic square}

Let $G$ be a complex
simple Lie group, $\fg$ its Lie algebra, and $G^{ad}$ the closed $G$-orbit
in
$\ppp
\fg$, the {\it adjoint variety} of $G$. 
Adjoint varieties are of current interest in algebraic geometry
 \cite{bea, hwang, lebrun}, it is conjectured that they are the only
complex contact manifolds with ample anticanonical bundle. 

In order to better understand the geometry of the adjoint varieties, one
could work infinitesimally. If one calculates the space of tangent
directions to lines passing through a point $x$ of $G^{ad}$, one obtains a
new variety $Y\subset \ppp T_xG^{ad}$. 
For example, in the case of $G=G_2$, $Y=v_3(\pp 1)\subset\pp 3$, the
twisted cubic curve. In other cases, to understand the geometry of $Y$
better, one can repeat the procedure. In the case of the remaining
exceptional groups, 
upon a second infinitesimalization  one arrives at the {\it
Severi varieties}, the projective planes over the composition algebras.
(These observations were
communicated to us by Y. Ye.)
The Severi varieties have been well studied, they arise in numerous
geometric
contexts. In particular, Zak \cite{zak} showed that they are the unique
extremal varieties for secant defects. They have the 
unusual  property that a generic hyperplane section of a Severi variety
is still homogeneous. Putting the resulting varieties into a chart we have:

\begin{center}\begin{tabular}{ccccc}
$v_2(Q^1)$ & $\PP (T\PP^2)$ & $G_\o(2,6)$ & $\BO\pp 2_0$ & \ hyperplane
section of Severi \\
$v_2(\PP^2)$ &
$\PP^2\times\PP^2$ & $G(2,6)$ & $\BO\pp 2$ &\ Severi \\ $G_\o(3,6)$ &
$G(3,6)$ &
$\SS_{12}$ & $E_7/P_7$ &\
lines through a point of $G^{ad}$\\ $F_4^{ad}$ & $E_6^{ad}$ & $E_7^{ad}$ &
$E_8^{ad}$ &\
$G^{ad}$
\end{tabular}\end{center}
where the notations are explained below. 

These varieties are homogeneous spaces of groups whose associated Lie
algebras are:

\begin{center}\begin{tabular}{cccc}
$\fso_3$ & $\fsl_3$ & $\fsp_6$ &
$\ff_4$\\ $\fsl_3$ & $\fsl_3\times \fsl_3$ &
$\fsl_6$ & $\fe_6$ \\ $\fsp_6$ & $\fsl_6$ & $\fso_{12}$ & $\fe_7$ 
\\ $\ff_4$ & $\fe_6$ & $\fe_7$ & $\fe_8$ \end{tabular}\end{center}
This chart is called {\em Freudenthal's magic square} of semi-simple Lie
algebras.

The magic square was constructed by Freudenthal and Tits as follows:
Let $\BA$ denote a complex composition algebra (i.e. the complexification
of $\BR$, $\BC$, the {\it quaternions} $\BH$ or the {\it octonions} $\BO$).
For a pair $(\AA,\BB)$ of
such composition algebras, the corresponding Lie algebra is $$\fg=Der\AA\op
(\AA_0\ot\JB_0) \op Der\JB, $$
where $\AA_0$ is the space of
imaginary elements, $\JB$ denotes the Jordan algebra of $3\times 3$
$\BB$-Hermitian
matrices, and $\JB_0$ is the subspace of $\JB$ consisting of traceless
matrices.
From Freudenthal's construction the symmetry in the chart appears to be as
miraculous as that $\fg$ is actually a Lie algebra. Vinberg gave a
construction where the symmetry is built in, see \cite{OV}.

The Severi varieties admit a common geometric interpretation as
$\BA\pp 2 \subset\ppp (\JA )$ and we showed in
\cite{lm0} that their hyperplane sections admit common geometric
interpetations as $G_Q(\BA^1, \BA^3)\subset\ppp (\JA_0)$, the
Grassmanian of $\BA^1$'s in $\BA^3$ isotropic for a quadratic form.
In particular, when one moves from left to right in the first two
rows, the varieties are naturally nested in each other. We show 
in \S 5,6 that the
same is true for the varieties above in the third and fourth rows,
in particular, we give a common geometric interpretation of the
varieties above in the third row as $G_w(\BA^3, \BA^6)$.

Moreover, as remarked above, as one moves from line to line
there are also natural inclusions of varieties (after fixing a
point and with the caveat that the inclusion of the first row
in the second is of a different nature).

The four Severi varieties share many common geometric properties:
their tangent spaces have a common geometric interpretation
as $\BA\oplus\BA$ (see \S  3),
and  the orbit structure in $\ppp (\JA)$,
the classification of  hyperplane sections,
and  the desingularizations
of the singular orbit closure
are all the same  (see \S 4,7,8). We show that these
extraordinary similarities also hold for the 
varieties above in the third and fourth rows.
Moreover, we show that the more complicated spaces of the third
and foruth rows can be understood in terms of a simple object,
the   discriminant of a cubic polynomial, as we now explain.

\subsection{The discriminant and generalizations}

  Consider  the twisted cubic curve $v_3(\pp 1)\subset \pp 3=
\ppp (S^3(\BC^2))$.
It is the space of cubic
polynomials having a triple root, and
its tangential variety, the
quartic hypersurface $\tau(v_3(\pp 1))\subset\pp 3$,
is the space of cubics having a multiple root. The equation $\Delta$
defining $\tau(v_3(\pp 1))$ is the classical {\em discriminant}
of a cubic polynomial
and is as follows: if
$P=p_0x^3+p_1x^2+p_2x+p_3$, then
$$\Delta (P) = 3(3p_0p_3-p_1p_2)^2+4(p_0p_2^3+p_1^3p_3)-4p_1^2p_2^2.$$

The ideal of $v_3(\pp 1)$ is generated by the second derivatives
of $\Delta$.
Write $W=\BC$ and write $\BC^4= V=\BC \oplus W\oplus W^*\oplus \BC^*$. 
Let $C(x)=x^3$. We rewrite the discriminant (changing scales) as
follows: for $w=( \a, r, s^*, \b^* )\in V$ let
$$Q(w)= (3\a\b^*-\frac{1}{2}\langle r,s^*\rangle )^2
+\frac{1}{3}(\b^*C(r^3)+\a C^*(s^{*3}))
-\frac{1}{6}\langle C^*(s^{*2}),C(r^2)\rangle .$$

We may describe  $v_3(\pp 1)$ as the image of the rational map:

$$\begin{array}{rcl}
\phi : \PP (\CC\op W) & \dashrightarrow & \PP(
V )=\PP(\CC\op W\op W^*\op\CC) \\
  ( z:w ) & \mapsto & (\frac{1}{6}z^3:z^2w:zC(w^2):\frac{1}{3}C(w^3)).
\end{array}$$

Letting $\fh= \fs\fl_2$,  
the adjoint variety   $G_2^{ad}\subset\ppp\fg_2$ is the image of the
rational map (see \S 5 and \cite{lm1}):

$$
\begin{array}{rcl}
\psi :
\ppp (\BC + V + \BC)&\ra & \ppp (\BC^*\oplus V^*\oplus (\BC\oplus \frak h)
\oplus V
\oplus \BC)\\
(u,A,v)&\mapsto &
(u^4, u^3A,   u^3v,  u^2Q(A,A,\cdot,\cdot ),
u^2vA - uQ(A, \cdot,\cdot, \cdot), u^2v^2 - Q(A)) .\end{array}
$$

\smallskip

Now, let $W  =\BC^{3m+3}$ be the vector space associated
to a Severi variety
$\BA\pp 2\subset\ppp W$.
$W$ is equipped with a cubic form $C$ (the determinant, see \cite{lm0}).
Let $Sp_6(\BA)$ (resp. $E(\BA)$) denote the groups appearing
in the third (resp fourth) row of the magic chart,
let $\fh = \fs\fp_6(\BA)$, and
  continue the notation $V= \BC \oplus W\oplus W^*\oplus \BC^*$
etc... We prove:

\smallskip

{\em 
The varieties
$G_w(\BA^3,\BA^6)\subset\ppp V$  are the images of the rational
mapping
$\phi$.
The quartic $Q$ is an $Sp_6(\BA)$-invariant
form  on $V=\BC^{6m+8}$. The hypersurface $Q=0$ is $\t (G_w(\BA^3,\BA^6))$
and the ideal of  $G_w(\BA^3,\BA^6)$ is
generated by the
second derivatives of $Q$. 
Moreover, the adjoint variety $E(\BA )^{ad}\subset\ppp (\fe (\BA))$ is
the image  of the rational map $\psi$.
}

\smallskip

The
 varieties $G_w(\BA^3,\BA^6)\subset\ppp V$ are also
Legendrian for a $Sp_6(\BA)$-invariant symplectic form $\Omega$
that generalizes the natural symplectic form on $S^3\BC^2$,
for which $v_3(\pp 1)$ is Legendrian, see \S 5.

\smallskip

We also show in \S 5 that in some sense   the quartic
$Q$ as well as Cayley's hyperdeterminant are  determined by
the classical discriminant.

The orbit structures for $G_w(\BA^3,\BA^6)\subset\ppp V$ are slighly more
complicated than for $v_3(\pp 1)$, as the first derivatives
of the quartic define an intermediate orbit closure which we
interpret as the locus of points on a family of secant lines,
see \S 5.3. We remark that   there are three   orbits
in the ambient space for the second row, four in the ambient space
for the third row and  five, not in $\ppp (\fe (\BA))$,
which contains an infinite number of orbits, but in the secant
variety of $E(\BA)^{ad}$.

\subsection{Notation}

$\BA$, $\BB$ denote complex composition algebras, i.e.
$\BA=\BA_{\BR}\ot_{\BR}\BC$ where $\BA_{\BR}= \BR, \BC, \BH$ or $\BO$ (the
four real division algebras). If $ a\in\BA$, $\overline a$ denotes its
conjugate as an element of $\BA$. $\BA\ot\BB$ denotes the tensor product
which has the algebra structure with multiplication
$(\a\ot b, a'\ot b')\mapsto aa'\ot bb'$
and conjugation $\overline{a\ot b}=\overline a\ot \overline b$. 

$\JA$ denotes
the space
of $\BA$-Hermitian matrices of order three, with coefficients in $\AA$:
$$\JA = \Biggl\{ \begin{pmatrix} r_1 & \overline{x_3} & \overline{x_2} \\
x_3 & r_2 & \overline{x_1} \\ x_2 & x_1 & r_3 \end{pmatrix}, \;\;
r_i\in\CC, \; x_j\in\AA \Biggl\}.$$ $\JA$ has the structure of a
{\it Jordan algebra} with the
multiplcation $A\circ B= \frac 12 (AB + BA)$ where $AB$ is the ususal
matrix multiplication. There is a well defined cubic form which we call the
{\it determinant} on $\JA$. 

We let $SO_3(\BA )\subset GL_{\BC}(\JA)$ denote the group of complex
linear transformations preserving the Jordan multiplication (the name is
motivated because the group can also be described as group preserving the
cubic form and the quadratric form $Q(A)=\trace (A^2)$). We have
respectively $SO_3(\BA ) = SO_3, SL_3, Sp_6, F_4$.

We let $SL_3(\BA )\subset GL_{\BC}(\JA)$ denote the group of complex linear
transformations preserving the determinant. Respectively $SL_3(\BA)= SL_3,
SL_3\times SL_3, SL_6, E_6$. 

${\cal Z}_2(\BA )$ denotes the space of {\it Zorn matrices}, $$
{\cal Z}_2(\BA )=
\left\{ \pmatrix x & A\\ B& y\endpmatrix \mid x,y\in \BC, A,B\in\JA\right\}
$$

It can be given the structure of an algebra, called a {\it Freudenthal
algebra},
 see \cite{kan}, \cite{fr}. $Sp_6(\BA)\subset Gl_{\BA}({\cal Z}_2(\BA ))$
respectively denotes the groups $Sp_6, SL_6, Spin_{12}, E_7$. It is the
group preserving the quartic discriminan on ${\cal Z}_2(\BA )$ (see
proposition \ref{disc}). \medskip

If $Y\subset\ppp T_xG/P$ is a subvariety, we let $\tilde Y\subset T(G/P)$
denote the corresponding distribution. 

If $X\subset\ppp V$, we let $\hat X\subset V$ denote the cone over $X$. 

When there is an orbit closure in $\PP V^*$ isomorphic to $X\subset\ppp V$,
we 
utilize
$X_*\subset\ppp V^*$ to denote this orbit closure.

\section{Freudenthal geometries}
\subsection{The magic square and the four geometries}

Freudenthal associates to each group in the square a set of preferred
homogeneous varieties ($k$ spaces for each group in the $k$-th row). These
spaces have the same incidence relations with the corresponding varieties
for the groups in the same row. He calls the geometries associated to
the groups of the rows respectively, {\em $2$-dimensional elliptic,
$2$-dimensional plane projective, $5$-dimensional symplectic} and {\it
metasymplectic}.
 The distinguished spaces are called respectively, spaces 
of {\it points, lines, planes and symplecta}. To avoid confusion, we will
use the terminology {\it F-points, F-planes} etc... 

\medskip
The spaces of elements are given by the following diagrams: 

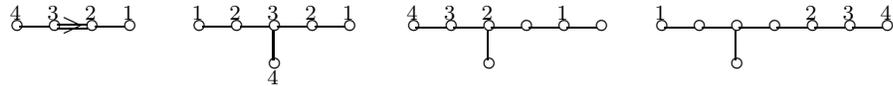
\begin{center}\begin{tabular}{cccc}
\setlength{\unitlength}{2.5mm}
\begin{picture}(10,3)(-2.5,-1)
\multiput(0,0)(2,0){4}{$\circ$}
\put(0.55,.35){\line(1,0){1.55}}
\put(4.55,.35){\line(1,0){1.55}}
\multiput(2.55,.25)(0,.2){2}{\line(1,0){1.55}} \put(2.75,0){$>$}
\put(0,.7){$\scriptstyle{4}$}
\put(2,.7){$\scriptstyle{3}$}
\put(4,.7){$\scriptstyle{2}$}
\put(6,.7){$\scriptstyle{1}$}
\end{picture} &
\setlength{\unitlength}{2.5mm}
\begin{picture}(10,3)(0,-1)
\multiput(0,0)(2,0){5}{$\circ$}
\multiput(0.55,.35)(2,0){4}{\line(1,0){1.55}}
\put(4.3,-1.4){\line(0,1){1.5}}
\put(4,-2){$\circ$}
\put(0,.7){$\scriptstyle{1}$}
\put(2,.7){$\scriptstyle{2}$}
\put(4,.7){$\scriptstyle{3}$}
\put(6,.7){$\scriptstyle{2}$}
\put(8,.7){$\scriptstyle{1}$}
\put(4,-2.7){$\scriptstyle{4}$}
\end{picture} &
\setlength{\unitlength}{2.5mm}
\begin{picture}(11,3)(0.8,-1)
\multiput(0,0)(2,0){6}{$\circ$}
\multiput(0.5,.3)(2,0){5}{\line(1,0){1.6}} \put(4,-2){$\circ$}
\put(4.3,-1.4){\line(0,1){1.5}}
\put(0,.7){$\scriptstyle{4}$}
\put(2,.7){$\scriptstyle{3}$}
\put(4,.7){$\scriptstyle{2}$}
\put(8,.7){$\scriptstyle{1}$}
\end{picture} &
\setlength{\unitlength}{2.5mm}
\begin{picture}(11,3)(0.8,-1)
\multiput(0,0)(2,0){7}{$\circ$}
\multiput(0.5,.3)(2,0){6}{\line(1,0){1.6}} \put(4,-2){$\circ$}
\put(4.3,-1.4){\line(0,1){1.5}}
\put(0,.7){$\scriptstyle{1}$}
\put(10,.7){$\scriptstyle{3}$}
\put(8,.7){$\scriptstyle{2}$}
\put(12,.7){$\scriptstyle{4}$}
\end{picture} \\
\end{tabular}\end{center}\medskip

\smallskip
Here a 1 denotes the space of F-points, 2 the space of F-lines, 3 the space
of F-planes, and 4 the space of F-symplecta in the metasymplectic
geometries. (E.g. the space of F-points for $E_6$ is $E_6/P_{1,6}$, where
we use the ordering of roots as in \cite{bou}, and $P_{1,6}$ is the
parabolic subgroup associated to the simple roots $\a_1$ and $\a_6$.)
Taking out the nodes
numbered $4$, we obtain the diagrams describing the three types of elements
in the $5$--dimensional symplectic geometries, and so on. 

\medskip
While Freudenthal was interested in the synthetic/axiomatic geometry 
of the spaces,
we are primarily interested in the spaces as subvarieties of a projective
space.
We have taken embeddings of the spaces to make the geometries as uniform as
possible.
Below are the spaces, which are in their minimal homogeneous embeddings
unless indicated. If $X\subset\ppp V$ is the minimal embedding,
$v_d(X)\subset\ppp S^dV$ indicates the $d$-th Veronese re-embedding. We use
standard nomenclature when there is one, and otherwise have continued the
labelling by parabolic.

\begin{center}\begin{tabular}{cccccccccc} $v_4(\PP^1)$ & $\BF^3_{1,2}$ &
$G_{\o}(2,6)$ & $\OO\PP^2_0$ & & & & & & \\ $v_2(\PP^2)$ &
$\PP^2\times\PP^2$ & $G(2,6)$ & 
$\OO\PP^2$ & & & $v_2(\PP^2)$ & $\PP^2\times\PP^2$ & $G(2,6)$ & $\OO\PP^2$
\\
$v_2(\PP^5)$ & $\FF^6_{1,5}$ & $G_o(2,12)$ & $E_7^{ad}$ & & & $G_{\o}(2,6)$
&
$\FF^6_{2,4}$ & $G_o(4,12)$ & $E_7/P_6$ \\ $v_2(\OO\PP^2_0)$ &
$E_6/P_{1,6}$ & $E_7/P_6$ & $E_8/P_1$ & & & $F_4/P_3$ & $E_6/P_{3,5}$ &
$E_7/P_4$ & $E_8/P_6$ \end{tabular} \end{center} 

\hspace*{3cm}{\bf F-points}\hspace*{6cm}{\bf F-lines} 

\begin{center}\begin{tabular}{cccccccccc} $G_{\o}(3,6)$ & $G(3,6)$ &
$\SS_{12}$ & $E_7/P_7$ & \hspace*{.5mm} & & & & & \\
$F_4/P_2$ & $E_6/P_{4}$ & $E_7/P_3$ & $E_8/P_7$ & \hspace*{.5mm} & &
$\hspace*{4mm}F_4^{ad} \hspace*{4mm}$ & $ \hspace*{4mm} E_6^{ad}
\hspace*{4mm}$ & $ \hspace*{4mm}E_7^{ad} \hspace*{4mm}$ & $
\hspace*{4mm}E_8^{ad}
\hspace*{4mm}$ \end{tabular} \end{center} 

\hspace*{3cm}{\bf F-planes}\hspace*{5.5cm}{\bf F-symplecta}

\medskip
Here we use the following notations: $G(k,l)$ denotes the Grassmanian of
$\BC^k$'s in $\BC^l$, $G_{\omega}(k,l)$ respectively $G_o(k,l)$ denotes the
Grassmanian of $\BC^k$'s in $\BC^l$ isotropic 
for a symplectic (resp. nondegenerate quadratic) form,
$\BF_{a,b}^c$ denotes the variety of
flags $\BC^a\subset\BC^b$ in a fixed $\BC^c$. $G^{ad}\subset\ppp\fg$
denotes the {\it adjoint variety} of $G$, the closed orbit in $\ppp \fg $. 

We will use the following notations: the F-points in respectively the
first, second and third rows, and columns corresponding to $\BA$,
will be denoted $\BA\pp 2_0, \BA\pp 2, Sp_6(\BA)^{ad}=G_w(\BA^1, \BA^6)$.
(In particular, $\OO\PP^2$ is the Cayley plane, see \cite{lm1}.) The
F-lines, resp. F-planes
in the third row will be denoted $G_w(\BA^2, \BA^6)$, resp. $G_w(\BA^3,
\BA^6)$. The notations are explained in \S 3.5. 

\begin{prop} Let $m=1,2,4,8$. The dimensions of the spaces of elements 
are as follows:
$$\begin{array}{|l|cccc|}\hline
 &\text{F-points} &\text{F-lines}
 &\text{F-planes}&\text{F-symplecta}\\ \hline
 \text{First row} &  2m-1  & & & \\
 \text{Second row} & 2m  &  2m  & &  \\
 \text{Third row} &   4m+1 &    5m+2  &  3m+3  &\\
\text{Fourth row} &  9m+6  &  11m+9  &  9m+11  &  6m+9  \\ \hline
\end{array}$$
\end{prop}

\medskip Freudenthal remarked that
for the fourth row, the usual duality
between elements of complementary dimension is lost already at the level of
their dimensions. We that vestiges of this duality remain. In particular,
giving the F-spaces \lq\lq F-dimensions\rq\rq\ $1,2,3,4$ respectively for
F-points, F-lines, F-planes, F-symplecta, if some geometric
element describes a space of dimension
$um+v$, the element of complementary F-dimension describes a space 
of dimension $vm+u$!

\subsection{The magic square for all $n$} 

One can define a magic square for all $n$, only one loses the fourth row
and column:

$$
\matrix
\fso_n & \fsl_n & \fsp_{2n}\\
\fsl_n& \fsl_n\times \fsl_n & \fsl_{2n}\\ \fsp_{2n} & \fsl_{2n} &
{\mathfrak spin}_{4n},\endmatrix 
$$
where now
$$
\fg= Der\BA\oplus (\BA_0\ot {\cal J}_n(\BB)_0)\oplus Der{\cal J}_n(\BB).
$$

We have  the following chart of F-points
$$
\matrix 
v_2(Q^{n-2}) & \BF^{n }_{1,n-1}&  G_{\omega}(2,2n)\\
v_2(\pp{n-1}) & Seg(\pp {n-1}\times\pp {n-1}) & G (2,2n)\\
v_2(\pp{2n-1})  &  \BF^{2n }_{1, 2n-1} & G_o(2,4n). \endmatrix
$$
The  analogue of F-planes, or perhaps better to say
F-hyperplanes, for the third row  are  the
minuscule varieties
$$
G_{\omega}(n,2n)\qquad  G(n,2n)\qquad \BS_{2n}.
$$

\noindent\rem If one allows a fifth column $\BA =0$ in the $n=3$ square
(making it into a rectangle), the additional column has Lie algebras
$0,0, \fso_3, \fg_2$ and thus one obtains all simple Lie algebras
from Freudenthal's magic.

\subsection{F-Schubert varieties and F-incidence}

Let $P$, $Q$ be parabolic subgroups of $G$. Consider the diagram
$$ \begin{array}{rcccl}
& &G & &\\
 & p\swarrow & & \searrow q& \\
X= G/P & & & & G/Q=Y\end{array}
$$
and define
$\Sigma^{Y, y_0}_X=\Sigma^{Y }_X=p(q\inv (y_0))= \{
x\in X\mid  x \ {\rm is\ incident\  to }\  y_0\}.$  
In the language of \cite{lm0}, $\Sigma^{Y, y_0}_X$ is the $(Y,X)$
Tits-transform of $y_0$.
We note that such a Schubert variety  furnishes a homogeneous vector bundle
over $Y$ by taking the fiber over $y\in Y$ to be the linear span of
the cone over  $\Sigma^{Y, y }_X$. 

The Freudenthal spaces distinguish certain Schubert
varieties which we will call {\em F-Schubert varieties}.   The  F-Schubert
varieties have uniform behavior as one changes $\BA$ and exhibit
similarities as
one changes the row. They play a role in understanding the geometries of
the F-varieties analogous to the role of classical Schubert
varieties for understanding Grassmanians. 

Let $m=1,2,4,8$.
In the case of the first row there is nothing to say. For the
second row, the variety of F-points incident to an F-line is
an $\BA\pp 1= Q^m$ and of course the variety of F-lines incident to an
F-point is
an $\BA\pp 1= Q^m$ as the two spaces are isomorphic.
 This symmetry
is broken with the third row as $\BA\pp 1=Q^m=G_o(1, m+2)$ generalizes
in two different ways, to $\BA\pp 2$ and to $G_o(2, m+4)$.
  \medskip

For the third row we have 
 $$ \begin{array}{|l|ccc|}\hline
  &\text{F-points} &\text{F-lines} &\text{F-planes} \\ \hline 
  \text{F-points} & &  \AA\PP^1 & \AA\PP^2\\
  \text{F-lines} & G_o(2, m+4) & & \AA\PP^2\\
  \text{F-planes} &  G_o(1,m+4) & \PP^1 &  \\ \hline
\end{array}$$

and for the fourth:

$$\begin{array}{|l|cccc|}\hline
 &\text{F-points} &\text{F-lines}
 &\text{F-planes}&\text{F-symplecta}\\ \hline
 \text{F-points} & &  \AA\PP^1 & \AA\PP^2 &G_{\o}(\AA^1,\AA^6) \\
 \text{F-lines} & G_o(3, m+6) & & \AA\PP^2&G_{\o}(\AA^2,\AA^6) \\
 \text{F-planes} & G_o(2,m+6) & \PP^2 & &G_{\o}(\AA^3,\AA^6)\\
\text{F-symplecta}&  G_o(1, m+6) &\PP^2 &\PP^1 & \\ \hline
\end{array}$$

The $\pp 2$'s corresponding to $\Sigma^{X_{F-lines}}_{X_{F-symplecta}}$
are embedded by the quadratic Veronese embedding.

\smallskip 
 
One can also study the incidence relations among elements of the same
space. Freudenthal  (\cite{fr2}, pp. 169-171) describes {\it uniform}
incidence
relations for the distinguished varieties of each row, which we
utilize in our study. See \S 6.

\section{Tangent spaces}

We describe the tangent bundles   in each of the Freudenthal
geometries $G/P$. Recall from \cite{lm0} that, as a module over 
a maximal semi-simple
subgroup of $P$, the decomposition of $T(G/P)$
into irreducible $H$-modules can be read off the root system of
$\fg$. Indeed, up to conjugation, the parabolic group $P$ is determined by
a set
of simple roots (a single root when $P$ is maximal), say $I$. Then the
irreducible components of the tangent bundle are, roughly speaking, in
correspondance with the possible coefficients of positive roots over the
simple
roots in $I$. We shall denote by $T_k$ the sum of the irreducible
components of
$T$ defined by coefficients over these simple roots with
the coefficients summing to $k$. 
 
\subsection{F-points}

 In the case of points, 
the tangent space has one or
two components, given by the following tables: 

\medskip\begin{center}
\begin{tabular}{cccclcccccc} 
${\bf T_1}$ & $\CC$ & $\CC\op\CC$ & $S^*\ot S^{\perp}/S$ & $\Delta$ & 
\hspace*{1cm} & ${\bf T_2}$ & * &
$\CC$ & $\CC^3$ & $\CC^7$\\ 
 & $\CC^2$ & $\CC^2\op\CC^2$ & $S^*\ot Q$ & $\Delta_+$ & \hspace*{1cm}
 & & * & * & * & * \\
&$\CC^4$ & $\CC^4\op\CC^4$ & $S^*\ot S^{\perp}/S$ & $\Delta_+$ &
\hspace*{1cm} & & $\CC$ &$\CC$ &$\CC$ &$\CC$ \\ 
& $\Delta$ & $\Delta\op\Delta$ & $S^*\ot \Delta$ & $\Delta_+$ &
\hspace*{1cm}
& & $\CC^7$ &$\CC^8$ &$\CC^{10}$ &$\CC^{14}$ \end{tabular} \end{center}

\begin{prop} Let $X=G/P$ be the space of
F-points in the Freudenthal geometry associated to the pair
of composition algebras $(\AA,\BB)$ and let $x\in X$.
Let $T_1\subset T_xX$ denote the smallest $P$-invariant sub-module. Let
$H$ be a maximal semi-simple subgroup of $P$.
Then $H$ is a spin group (or product of spin groups) and
$T_1$ is a spin representation (or product of such). Moreover,
$$T_1\simeq \AA\ot\BB.$$
\end{prop}

\begin{proof}
Let $\Delta^k=\Delta^{B_k}$, $\Delta_+^k=\Delta_+^{D_k}$,
$\Delta_-^k=\Delta_-^{D_k}$. We rewrite the table for $T_1$ as

\medskip\begin{center}
\begin{tabular}{ccclccccc} 
$\Delta^0$ & $\Delta^0_+\op\Delta^0_-$ & $\Delta^1_+\ot\Delta^1$ &
$\Delta^3$ &
\hspace*{1cm} &   &   &   & \\ 
$\Delta^1$ & $\Delta^1_+\op\Delta^1_-$ & $\Delta^1_+\ot\Delta^2_+$
 & $\Delta^5_+$ &
\hspace*{1cm} &   &   &   &   \\
$\Delta^2$ & $\Delta^3_+\op\Delta^3_-$ & $\Delta^1_+\ot\Delta^4_+$
 & $\Delta^6_+$ &\hspace*{1cm} &   &   &   &  \\ 
$\Delta^3$ & $\Delta^3_+\op\Delta^3_-$ & $\Delta^1_+\ot\Delta^5_+$
 & $\Delta^7_+$& \hspace*{1cm} &   &   &   & 
\end{tabular} \end{center}

\medskip The fact that several spin representations of small 
dimensions have natural realizations given by composition algebras
can be found in \cite{harvey}, from which all cases except
for  $\BO\ot\BO$ 
can be deduced.  The  case of  
$\BO\ot\BO$ follows from   proposition 3.3 below, which
gives a general relation  between Clifford algebras and composition
algebras.\qed

\medskip
Note that in the magic chart for $n>3$, the tangent space to points
in the third row does not have an analogous interpretation.\end{proof}

\medskip 
The space $T_2$ for points does not have a very regular behavior.
Note however that for the first and last lines $T_2$
has the interpretation of $\BA_0 \op\BB_0$. While $T_2$ does not
behave well for points, we have the following proposition:

\begin{prop} Notations as above. For each Lie algebra $\fg$
in   Freudenthal's
magic square, 
there is a parabolic subgroup $\fp$ of $\fg$  such that
the quotient $\fg /\fp$ decomposes into the sum of $$T_1\simeq \AA\ot\BB
\quad and\quad T_2\simeq \AA_0\op\BB_0,$$
where $\fg$ is associated to the pair
of composition algebras $(\AA,\BB)$.

$H$ is a spin group (or product of spin groups),
$T_1$ is a spin representation (or product of such)
and $T_2$ is a vector representation. 
 \end{prop}

\medskip \noindent For the first and last lines of the square,  
the corresponding $G/P$'s are the spaces of F-points, while for the third
line 
they are the spaces of F-lines. For the second line they are the spaces of
incident pairs of F-points and F-lines. Note in particular that the square
below, formed by these $G/P$ is perfectly symmetric, although the geometric
interpretations are not. \medskip

\begin{center}\begin{tabular}{cccc}
$v_4(\PP^1)$ & $\FF_{1,2}$ & $G_\o(2,6)$ & $\OO\PP^2_0$ \\ $\FF_{1,2}$ &
$\FF_{1,2}\times\FF_{1,2}$ & $\FF_{2,4}$ & $E_6/P_{1,6}$ \\ $G_\o(2,6)$ &
$\FF_{2,4}$ & $G_o(4,12)$ & $E_7/P_6$\\ $\OO\PP^2_0$ & $E_6/P_{1,6}$ &
$E_7/P_6$
& $E_8/P_1$ \end{tabular}\end{center}\medskip

\begin{prop} Let $\BA, \BB$ be complex composition algebras, other
  than the complexification of $\RR$.
Let $\AA_0\oplus \BB_0$ be endowed with the quadratic form
$Q(a+b )= a\overline a  + b\overline b= -a^2- b^2$.
Then there is a natural diagram of maps  of algebras
$$\begin{array}{ccc}
Cl(\AA_0\oplus \BB_0, Q) & \ra &  End (\BA\ot\BB \op \BA\ot \BB) \\
\cup & & \cup \\
Cl(\AA_0\oplus \BB_0, Q)^{even} & \ra &  End (\BA\ot\BB)\ot End
(\BA\ot\BB).  
\end{array}
$$
\end{prop}

The inclusion on the right is the \lq\lq diagonal\rq\rq\ inclusion. When
$\BA=\BB=\OO$ the dimensions of $Cl(\OO_0\oplus \OO_0, Q)^{even}$ and  
$ End (\OO\ot\OO)\ot End (\OO\ot\OO)$ coincide, showing that the 
two half-spin representations of $Spin_{14}$ have natural realizations 
on $\OO\ot\OO$. \medskip

\begin{proof}
By the fundamental lemma of Clifford algebras, see \cite{harvey}, 
we have to construct a map
$$
\phi :\AA_0\oplus \BB_0\ra End (\BA\ot\BB \op \BA\ot \BB)
$$
such that $\phi (a+b)^2= Q(a+b)Id_{\BA\ot\BB \op \BA\ot \BB}$,
as then there exists a unique extension to a map
$\tilde\phi : Cl(\AA_0\oplus \BB_0, Q)\ra End (\BA\ot\BB \op \BA\ot \BB)$.
Consider 
$$
\phi (a+b) (\a\ot\b , \c\ot\d )=
(ia\c\ot\overline\d + \c\ot\overline \d b,
ia\a\ot\overline\b + \a\ot b\overline\b).
$$
 A short calculation shows that $\phi$ has the required property.
Moreover, since $Cl^{even}$ is generated by the products of even numbers 
of vectors in $\AA_0\oplus \BB_0$, the diagram follows. 
\end{proof}\qed

\subsection{F-lines} 

The components of the tangent spaces for F-lines are as follows:

\begin{center}\begin{tabular}{ccccccccccc} 
${\bf T_1}$ & $\CC^2$ & $\CC^2\op\CC^2$ & $\CC^2\ot\CC^4$ & $\Delta_+$
& \hspace*{1cm} & ${\bf T_3}$ &
* & * & * & *\\
& $\CC^2\ot\CC^2$ & $\CC^2\ot\CC^2\op\CC^2\ot\CC^2$ & $\CC^4\ot\CC^4$ &
$\CC^2\ot\Delta_+$
& \hspace*{1cm} & & * & *& * & *\\
& $\CC^2\ot\CC^3$ & $\CC^2\ot\CC^3\op\CC^2\ot\CC^3$ &
$\CC^2\ot\CC^3\ot\CC^4$
& $\CC^3\ot\Delta_+$ & \hspace*{1cm} & & $\CC^2$ & $\CC^4$ & $\CC^8$ &
$\CC^{16}$ \\ 
& & & & & & & & & & \\
${\bf T_2}$ & * &* &* &* & \hspace*{1cm} & ${\bf T_4}$ &*
&* &* &* \\ 
& $\CC^3$ & $\CC^4$ & $\CC^6$ & $\CC^{10}$ & & &* &* &* &* \\ &$\CC^9$ &
$\CC^{12}$ & $\CC^{18}$ & $\CC^{30}$ & \hspace*{1cm} & & $\CC^3$ & $\CC^3$
&
$\CC^3$ & $\CC^3$ \end{tabular} \end{center}

Let $Y_{\BA}$      respectively denote
$\emptyset$, $\pp 1\sqcup \pp 1$, $\pp 1\times\pp 3$, $\BS_5$ and
let $H_{\BA}$  respectively denote
$SL_2$, $SL_2\times SL_2$,
$SL_2\times SL_4$ and $Spin_{10}$.
Note that $Y_{\BA}\subset\ppp (\BA\oplus \BA)$ and if we
give $\BA\oplus\BA$ coordinates $(u,v)$ then
$I_2(Y_{\BA})=\{ u\overline u, v\overline v, u\overline v\}$, having
respectively
$3,4,6$ and $10$ generators. Examining the spaces above, we obtain:

\begin{prop} Let $X=X_{F-lines}^{p,\BA}=G/P$ denote the space of F-lines
in the $p$-th row whose  composition algebra is
$\BA$. With the same notations as above,
$$
\begin{array}{rcl}
T_1 &=&\BC^{p-1}\ot \BA^2\\
Y_1 &=& Seg (\pp{p-2}\times Y_{\BA})\\
H &=&SL_{p-1} \times H_{\BA}\\
 T_2&=&\La 2\BC^{p-1}\ot I_2(Y_{\BA}).\end{array}
$$
\end{prop}

\subsection{F-planes}

\begin{prop} Let $X=X_{F-planes}^{p,\BA}=G/P$ denote the space of F-planes
in the $p$-th row  whose composition algebra is  $\BA$. With the same
notations as above, 
$$
\begin{array}{rcl}
T_1 &=&\BC^{p-2}\ot \JA\\
Y_1 &=& Seg (\pp{p-3}\times \BA\pp 2)\\
H &=&SL_{p-2} \times SL_3{\BA}\end{array}
$$
For $p=3$, $T=T_1$ and for $p=4$, 
$T_2={\cal J}_3(\AA)^*$ and
$T_3=\CC^2.$ 
\end{prop}
 
Unlike the case of F-points,
the structure of the tangent space for F-hyperplanes has a similar
interpretation for all
$n$. As with the $n=3$ chart, the tangent directions to
lines through a point is the space of F-points of the second row.
Let $\BA\pp n$ respectively 
denote $  v_2(\pp{n-1}), Seg(\pp{n-1}\times
\pp{n-1}), G(2,2n)$ and let $SL_n(\BA)$ respectively
denote $SL_n, SL_n\times SL_n, SL_{2n}$: 

\begin{prop}  Let $X=X_{F-planes}^{ \BA}=G/P$ denote the space of
F-hyperplanes
    in the $3$-rd row of the generalized chart whose composition algebra is

$\BA$. With the same notations as above, $T_2=0$ and 
$$
\begin{array}{rcl}
T_1 &=&  {\cal J}_n(\BA) \\
Y_1  &=&  \BA\pp {n-1}\\
H &= &   SL_n({\BA}).\end{array}
$$
\end{prop}

\subsection{F-symplecta}

The spaces of F-symplecta are  the adjoint varieties 
of the exceptional groups other than $G_2$.  

\begin{prop} Let $X=X_{F-symplecta}^{\BA}=G/P$ denote the space of
F-symplecta  whose composition algebra is $\BA$. With the same
notations as above, 
$$
\begin{array}{rcl}
T_1 &=&{\cal Z}_2(\AA)\\
T_2&=&\BC \\
Y_1 &=& G_w(\BA^3, \BA^6)\\
H &= & Sp_6({\BA})\end{array}
$$
\end{prop}
 
\subsection{Interpretations as Grassmanians}

Let $X=G/P\subset\ppp V$ be a homogeneous variety  with $P$ maximal.
When $G$ is a classical group, $X$ can be characterized as a family 
of $k$-planes in some natural representation of $G$. We investigate 
the existence of similar characterizations in the exceptional cases,
in terms of composition algebras (for a different kind of such 
characterizations, see \cite{lm1}, Corollary 7.8).
For the varieties of F-points in the second row of the magic square, 
$X=G(\BA^1, \BA^3)=\BA\pp 2$ and
in \cite{lm0} we gave the interpretation of the varieties
of F-points in the first row as $G_Q(\BA^1,\BA^3)=
\BA\pp 2_0$. (Being null for the cubic and the trace is equivalent to 
being null for the
cubic and  quadratic forms $Q(x)=\ttrace (x\circ x)$
where $\circ$ is the Jordan multplication in $\JA$.)

Rozenfeld announces  (\cite{rose}, theorem 7.22)  a unified geometric
interpretation of certain varieties in  the chart, which he calls 
elliptic planes over $\BA\ot\BB$. Also, before Rozenfeld's work was 
available in English, E. Weinstein \cite{w} had conjectured
that there should be a unified interpretation of 
varieties in the chart as some type of Grassmanian over $\BA\ot\BB$.  

As a first step towards Weinstein's conjecture, we felt that just as the
tangent space to an ordinary Grassmanian has an interpretation as $T_EG(k,
V) =E^*\ot V/E$, for there to be a unified interpretation of the chart,
there should be a unified interpretation of tangent spaces, and this
infinitesimal problem is solved above. In what follows we suggest global
interpretations based on our infinitesimal calculations. 

\medskip
  We begin with the F-varieties of the third row. 
If $P\in\JA$, its {\it comatrix} is defined by
$$\com(P)=P^2-(\tr P)P+\frac{1}{2} ((\tr P)^2-\tr P^2)I,$$
and characterized by the identity $\com(P)P=det(P)I$. Thus
the linear form $P\mapsto\trace (\com(P)P)$ is a polarization of the 
determinant. The varieties of F-planes in the third
row are the image of the rational map $\phi$ 
described in \S 1.2.  On an affine open subset 
we may write 
$$\phi (1,P)=(1,P,\com (P), \det (P)).$$
Note in particular we recover the natural identification    
 $T_xX^{3,\BA}_{F-planes}\simeq \JA$ from the map $\phi$ alone. Moreover,
if
$(I_3,P)$ is  interpreted as a matrix of three
row  vectors in $\AA^6$, the map
$\phi$ is  the usual Plucker map. The 
condition that $P\in\JA$ can be interpreted as the fact that the three
vectors defined by the matrix $(I,P)$ are orthogonal with respect to
the Hermitian symplectic two-form 
$w(x,y)={}^txA\overline{y}$,
where $A=\bigl(\begin{smallmatrix} 0 & I \\ -I & 0 
\end{smallmatrix} \bigl)$. It is therefore 
natural to see   $X^{3,\BA}_{F-planes}$ as a 
of Grassmannian of symplectic three-planes in $\AA^6$. This 
motivates our notation 
$$X^{3,\BA}_{F-planes}=G_{\o}(\AA^3,\AA^6).$$

Similarly, consider a  matrix of two vectors in $\AA^6$, 
of the form $(I_2,R,S)$, where $R$ and $S$ are two matrices of order
two. These two vectors are orthogonal with respect to
the Hermitian symplectic two-form $w$ if and only if $S$ is
Hermitian. The space of F-lines can be interpreted as a
Grassmannian $G_{\om}(\AA^2,\AA^6)$ of symplectic two-planes in 
$\AA^6$, and its tangent space, as expected, decomposes into
$T_1=\HH\ot\AA$ and $T_2={\cal J}_2(\AA) $. 

The space of F-points also has an interpretation as
$G_{\om}(\AA^1,\AA^6)$, its tangent space decomposing into
$T_1=\HH\ot\AA$  and $T_2=\CC$.

  In summary: 

\begin{prop} The spaces of F-points, F-lines and F-planes for the 
third line of the magic square, have a natural interpretation as 
symplectic Grassmannians $G_{\om}(\AA^k,\AA^6)$, with $k=1,2,3$ 
respectively.   Their tangent spaces 
are respectively  $\HH\ot\AA\op\CC$, $\HH\ot\AA\op {\cal J}_2(\AA)$
and $\JA$. \end{prop} 

\medskip
\noindent{\bf Problem}. {\em Find a unified interpretation of the
F-varieties of the fourth row.}

\medskip
Consider the quadratic form on $(\BA\ot\BB)^3$, with values in
$\AA\ot\BB$,  given for
$x=(x_1,x_2,x_3)$  by $Q (x )=x_1\overline x_3 + x_2\overline x_2$,
and consider the space
$$
G_Q(\BA\ot \BB, (\BA\ot\BB )^3) :=
\{ x\in (\BA\ot\BB )^3 \mid Q(x)=0\}.
$$
Then $T_xG_Q(\BA\ot \BB, (\BA\ot\BB )^3) =\{  y \mid Q(x,y) =0\}.$
If $x=(1, 0,0)$, 
we need $y_2\in Im (\BA\ot \BB)\simeq \BA_0\oplus \BB_0$
and there is no restriction on $y_3\in\BA\ot\BB$. This suggests
that the varieties in proposition 4.3 
({\em not} the varieties of F-points) admit a common interpretation
as $G_Q(\BA\ot \BB, (\BA\ot\BB )^3)$.

\medskip
Regarding F-varieties for $n>3$, we have the following proposition:

\begin{prop} Let $n>3$. We have the following interpretations:

F-points of the first row: $\BA\pp{n-1}_0$

F-points of the second row: $\BA\pp{n-1} $

F-hyperplanes of the third row: $G_w(\BA^n,\BA^{2n})$.

 In particular,
we have the identification 
$\BS_{2n}=G_w(\BH^n,\BH^{2n})$.
\end{prop}

\begin{proof} Note that
$SL_n(\BH )= SL_{2n}$ to obtain that the points of the
second row are indeed $\BA\pp {n-1}$ and the first are
$\BA\pp{n-1}_0$ as in the four by four case. For the third row, one
uses the same argument as above, only note that the corresponding Plucker
type mapping is of degree $n$.\end{proof}\qed

 \section{Folding and hyperplane sections of Severi varieties}

\subsection{Severi varieties}

We have little new to say about the F-points (or F-lines) of the second
row, otherwise known as the {\it Severi varieties}
$\BA\pp 2\subset\ppp\JA$ which may be described as the projectivization
of the rank one elements of $\JA$. Their secant varieties
$\sigma (\BA\pp 2 )$ are the rank at most two elements, i.e.,
the hypersurface $det =0$. Throughout this section we let
$m=1,2,4,8$.

  We record the following known proposition:

\begin{prop} $SL_3(\BA)$ has three orbits on $\BP\JA$, namely
$\BA\pp 2$, $\sigma (\BA\pp 2)\backslash \BA\pp 2 $,  and the 
open orbit $\PP\JA\backslash \sigma (\BA\pp 2)$, 
 which respectively correspond to
the matrices of rank one, two and three. \end{prop}

The unirulings of $\BA\pp 2$
are described in \cite{lm0}. They are all $SL_3(\BA )$-homogeneous,
i.e.,  given by Tits
transforms.  Here we give several descriptions of  the
 ruling  of  $\s (\BA\pp 2)$    by $\pp{m+1}$'s:

 The   rulings of $\s (\BA\pp 2)$ were implicitly described by Zak 
as follows: for $p\in \s(\BA\pp 2)\backslash \BA\pp 2$, let
$$
\Sigma_p=\overline { \{ x\in \BA\pp 2 \mid \exists y\in\BA\pp 2
\ {\rm such}\ {\rm that}\ p\in\pp 1_{xy}\} }
$$
the {\it entry locus of $p$}. Then Zak
\cite{zak} shows that $\Sigma_p$
is a quadric hypersurface in a $\pp{m+1}$, i.e.,
an  $\BA\pp 1$,  and that  $\s (\BA\pp 2)$
is therefore  ruled by these $\pp{m+1}$'s.  

 Another way to
view this  ruling is as follows: Let $\BA\pp 2_*\subset\ppp
\JA^*$ denote the closed orbit in the dual projective space. Then
$\s (\BA\pp 2)^*= \BA\pp 2_*$, i.e. the dual of the secant variety 
of $\BA\pp 2$ is the   the closed orbit in the dual
projective space (and $\s (\BA\pp 2)= (\BA\pp 2_*)^*$
by the reflexivity theorem). Let $N^*$ denote the conormal bundle to 
$\BA\pp 2_*$. Given $H\in \BA\pp 2_*$, 
$$\BP N^*_H=\{ p\in \s (\BA\pp 2) \mid
\tilde T_p\s (\BA\pp 2)\subset H\}$$ is the corresponding
$\pp{m+1}$ for any $p\in\BP N^*_H\backslash \BA\pp 1\subset
\ppp\JA $.

 The  rulings may also be seen from
Freudenthal's perspective:   $\BA\pp 2$
is the space of F-points and $ \BA\pp 2_*$ is the space of F-lines.
The F-Schubert variety  
of a $p\in \BA\pp 2_*$ is an $\BA\pp 1=Q^m
\subset\BA\pp 2$ and this $\BA\pp 1$ is the variety describe above.

\medskip
The Severi varieties were constructed by Zak using a
degree two map defined by the quadrics vanishing
on $Y_{\BA}\subset\PP(\AA\op\AA)=\pp{n-1}\subset\pp n$, see \cite{zak}. 
This construction can be generalized to construct the varieties
of F-lines in the third and fourth rows as well, using
$Y=Seg (\pp 1\times Y_{\BA} )$ for the third row and 
$Y=Seg (\pp 2\times Y_{\BA})$ for the fourth row. See
\cite{lmafter}. 
 
\subsection{Geometric folding}

To deduce the first line of the magic square from the second
one  we use the {\em folding} of a root system. Consider some Dynkin 
diagram with a two-fold symmetry $\th$, and let $\fg$ be the corresponding 
simple Lie algebra. 
If we choose a system of Chevalley generators for
$\fg$, there is a uniquely defined algebra involution of $\fg$ 
inducing the automorphism $\th$ of the simple roots \cite{OV}. 
Hence a decomposition $\fg=\ff\op W$ into eigenspaces, where 
$\ff$ is a Lie subalgebra, and $W$ has a natural $\ff$-module 
structure. A case-by-case examination then gives:

\begin{prop} There is a commutative (in general non-associative) 
$\ff$-equivariant 
multiplication on $W$,  and $V=\BC\op W$ is a simple $\fg$-module.
\end{prop} 

In most cases $V$ inherits an algebra
structure  
for which $F\subset GL(W)\subset GL(V)$ is the group preserving 
the structure.   
Here is a chart summarizing the representations arising from folding: 

\setlength{\unitlength}{3mm}
\begin{picture}(25,19)(2,-12.5)

\multiput(9,0)(2,0){3}{$\circ$}
\put(9,0){$\bullet$}
\put(9,-6){$\bullet$}
\multiput(9.45,.3)(2,0){2}{\line(1,0){1.65}}
\multiput(14.7,1)(2,1){1}{$\circ$}
\multiput(14.7,-1)(2,-1){1}{$\circ$}
\put(13.5,.5){\line(2,1){1.3}}
\put(13.55,0.1){\line(2,-1){1.3}}
\put(9,2.8){$D_{n+1}$}
\multiput(9,-6)(2,0){4}{$\circ$}
\multiput(9.45,-5.7)(2,0){2}{\line(1,0){1.65}}
\multiput(13.45,-5.8)(0,0.2){2}{\line(1,0){1.65}} \put(13.7,-6){$>$}
\put(9,-4){$B_{n}$}

\multiput(20,0)(2,0){2}{$\circ$}
\put(20.45,.3){\line(1,0){1.65}}
\multiput(23.7,1)(2,1){2}{$\circ$}
\put(25.7,2){$\bullet$}
\put(26,-6){$\bullet$}
\multiput(23.7,-1)(2,-1){2}{$\circ$}
\put(22.5,.5){\line(2,1){1.3}}
\put(24.3,1.45){\line(2,1){1.4}}
\put(22.55,0.1){\line(2,-1){1.3}} 
\put(24.3,-0.9){\line(2,-1){1.4}} 
\put(20,2.8){$E_6$}
\multiput(20,-6)(2,0){4}{$\circ$}
\multiput(20.45,-5.7)(2,0){1}{\line(1,0){1.65}}
\multiput(24.45,-5.7)(2,0){1}{\line(1,0){1.65}} 
\multiput(22.4,-5.8)(0,0.2){2}{\line(1,0){1.7}} \put(22.7,-6){$>$}
\put(20,-4){$F_4$}

\put(31,2.8){$A_{2n-1}$}
\multiput(31,0)(2,1){4}{$\circ$}
\put(33,1){$\bullet$}
\put(33,-6){$\bullet$}
\multiput(31,0)(2,-1){4}{$\circ$}
\multiput(31.6,.45)(2,1){3}{\line(2,1){1.4}}
\multiput(31.6,.15)(2,-1){3}{\line(2,-1){1.4}}
\multiput(31,-6)(2,0){4}{$\circ$}
\multiput(33.45,-5.7)(2,0){2}{\line(1,0){1.65}}
\multiput(31.45,-5.8)(0,0.2){2}{\line(1,0){1.65}} \put(31.7,-6){$>$}
\put(31,-4){$C_n$}

\put(42,2.8){$A_n\times A_n$}
\multiput(42,.7)(2,0){4}{$\circ$}
\multiput(42,-.7)(2,0){4}{$\circ$}
\multiput(42.45,1)(2,0){3}{\line(1,0){1.65}}
\multiput(42.45,-.4)(2,0){3}{\line(1,0){1.65}}
\put(42,.7){$\bullet$}
\put(48,-.7){$\bullet$}
\multiput(42,-6)(2,0){4}{$\circ$}
\multiput(42.45,-5.7)(2,0){3}{\line(1,0){1.65}}
\put(42,-6){$\bullet$}
\put(48,-6){$\bullet$}
\put(42,-4){$A_n$}

\put(5,-10){$W$}
\put(5,-12){$V$}
\put(10,-10){$\CC^{2n+1}$}
\put(10,-12){$\CC^{2n+2}$}
\put(21,-10){$\JA_0$}
\put(21,-12){$\JA$}
\put(30,-10){$\Lambda^{\langle 2\rangle}\CC^{2n}={\cal J}_n(\BH)_0$}
\put(31,-12){$\Lambda^2\CC^{2n}={\cal J}_n(\BH)$}
\put(44,-10){$\fsl_n$}
\put(43,-12){$M_n(\CC)$}
\end{picture}

\medskip\noindent\begin{rema} 
It follows from the results in \cite{akh} that if 
$X\subset\BP V$ is a homogeneous variety (under a semi-simple group) 
such that a generic hyperplane section of $X$ is still homogeneous, and if
$X\subset\BP V$ is not $\pp m\subseteq\pp m$ or $Q^m\subset\pp{m+1}$
(the two self-reproducing cases), then it must be the 
variety of F-points in the second row of a magic chart and $X\cap H\subset
H$ is the corresponding variety of F-points in the first row.

Note that there is a slight anomaly in that the chart does not
exactly correspond to geometric folding. However the
exception, $v_2(Q)$ can by seen as a special case of a different
phenomenon:
{\it any} homogeneous variety $X=G/P\subset\ppp V$ can be realized
as $v_2(\ppp V)\cap \ppp (S\up 2 V)$. In the case of $v_2(Q)$,
$\ppp (S\up 2 V)$ happens to be a generic hyperplane.
Note also that the $Q^{2m}\cap H$ section
is accounted for by geometric folding but the $Q^{2n+1}\cap H$ section
is not.

\end{rema}

\medskip Let $H_F\subset H_G$ denote the corresponding maximal semi-simple
subgroups of the isotropy groups fixing a point of $x\in X\cap H\subset X$.
 In the case of $X=\OO\PP^2$,   $H_G= Spin_{10}$   acts 
irreducibly on $TX=\Delta_+$, a sixteen-dimensional half-spin 
representation.  As an $H_F= Spin_7$-module    $T_x
 \BO\pp 2_0=\Delta\op\CC^7$, the sum of the 
spin and  the vector representations. 

This decomposition is interesting because spin 
representations of spinor groups 
usually decompose into sums of spin representations
 when restricted to smaller spinor
groups. The appearance of the vector representation may be understood
in terms of  the    
triality  
automorphism of $Spin_8$. We consider
$Spin_7\subset Spin_8$. The relevant embedding   $Spin_8\subset
Spin_{10}$ is such that, because of triality, $\Delta^{D_5}_+$ decomposes
as
$\Delta^{D_4}_+\oplus V^{D_4}$ as a $D_4$-module. When one
restricts further to $Spin_7\subset Spin_8$, $V^{D_4}=
V^{B_3}\oplus\BC$ and $\Delta^{D_4}_+= \Delta^{B_3}$.

 One can see the situation pictorially by considering
the fold of the Dynkin diagram 
of $E_6$ into the diagram of $F_4$: 

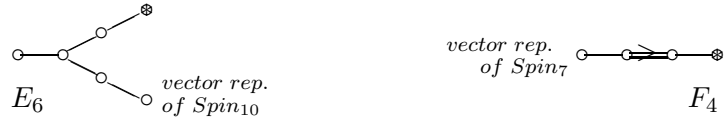
\setlength{\unitlength}{3mm}
\begin{picture}(15,9)(-10,-5)

\multiput(0,0)(2,0){2}{$\circ$}
\put(.45,.3){\line(1,0){1.65}}
\multiput(3.7,1)(2,1){2}{$\circ$}
\put(5.7,2){$*$}
\multiput(3.7,-1)(2,-1){2}{$\circ$}
\put(2.5,.5){\line(2,1){1.3}}
\put(4.3,1.45){\line(2,1){1.4}}
\put(2.55,0.1){\line(2,-1){1.3}} 
\put(4.3,-0.9){\line(2,-1){1.4}}
\put(6.7,-1.2){${\scriptstyle vector\;rep.}$}
\put(6.7,-2.2){${\scriptstyle of\;Spin_{10}}$}
\put(0,-2){$E_6$}
 
\multiput(25,0)(2,0){4}{$\circ$}
\put(31,0){$*$}
\multiput(25.45,.3)(2,0){1}{\line(1,0){1.65}}
\multiput(29.45,.3)(2,0){1}{\line(1,0){1.65}} 
\multiput(27.4,.2)(0,0.2){2}{\line(1,0){1.7}} 
\put(27.7,.1){$>$}
\put(30,-2){$F_4$}
\put(19.3,.5){${\scriptstyle vector\;rep.}$}
\put(20.8,-.5){${\scriptstyle of\;Spin_7}$}

\end{picture}

\section{F-planes in the third row}

\subsection{Constructions}

A special case of the
minuscule algorithm in \cite{lm1} constructs
$G_w(\BA^3, \BA^6)\subset\ppp ({\cal Z}_2(\BA))$ from $\BA\pp 2$ via a
degree
three mapping,  as well as constructing
$\fsp_6(\BA)$ from $\fsl_3(\BA)$ as
$\fsp_6(\BA)=\BC^*\oplus \JA^*\oplus
(\fsl_3\BA + \BC) \oplus\JA\oplus \BC$ . The construction also produces the
increasing 
filtration of
${\cal Z}_2(\BA)$ as a $U(\fsl_3\BA)$-module, namely
$$
{\cal Z}_2(\BA)=\BC\oplus \JA \oplus \JA^*\oplus\BC.$$
The action of $\fsp_6(\BA)$ can also  be described in
terms of creation and annhilation, see \cite{lm1}.

\medskip\noindent
\begin{rema}The F-planes in the fourth row can be constructed by a  
mapping
defined by polynomials vanishing on $Y_1=\pp 1\times\BA\pp 2$, $Y_2=
\BA\pp 2_*$ and their auxilliary varieties.
 See \cite{lm1}, \cite{lmafter} for details.\end{rema}

\medskip
If one allows the fifth algebra $\BA=\underline 0$, so that 
$$
{\cal J}_3(\underline 0)=
\{
\pmatrix r_1 & & \\ & r_2 & \\
 & & r_3\endpmatrix \mid r_j\in \BC\}, \qquad
 {\cal Z}_2(\underline 0)
=\left\{\pmatrix a & X
\\  Y & b\endpmatrix
\mid   a,b\in\BC,\ X,Y\in  {\cal J}_3(\underline 0) \right\}, 
$$
  the above construction
works equally well, except for the unfortunate notations
$Sl_3(\underline 0)= \BC^* + \BC^*$,
$\fsp_6(\underline 0)= \fsl_3 + \fsl_3 + \fsl_3$, 
$\underline 0\pp 2= \pp 0\sqcup \pp 0\sqcup\pp 0$, $G_w(\underline 0^3,
\underline 0^6)= Seg(\pp 1\times\pp 1\times \pp 1)$.

\medskip
Yet another space will play an important role: let
$\Delta_{\cal J}\subset {\cal J}_3(\underline 0)$ denote the
homotheties and let
$\Delta_{\cal Z}\subset {\cal Z}_2(\underline 0)$ denote the
subspace induced by $\Delta_{\cal J}$, i.e., the  subspace where
$X$ and $Y$ are homotheties. 

Here is  an alternate construction of
 $\fsp_6(\BA)$ and ${\cal Z}_2(\BA)$  that makes
no reference to composition algebras, one only  uses
the existence of an invariant cubic polynomial:

\begin{theo} (geometric version) Let $Z=H/Q\subset\ppp W$ be a homogeneous
variety with $H$ simple having the properties that
closure of  the largest $H$-orbit  in $\ppp W$ is a cubic hypersurface
and that $I_2(Z)$ is an irreducible $H$-module.

 Then $\fg=W^*\op  (\fh +\BC ) \op W$ is a 
simple Lie algebra and $V=\CC\op W \op W^* \op\CC$ has a natural
structure of simple $\fg$-module.  
Moreover, if $X\subset\ppp V$ denotes the closed $G$-orbit, then the
space of $\pp 1$'s in $X$ through a point $x$ is isomorphic to
$Z$.\end{theo}

The significance of this theorem is due to the set of varieties satisfying
its hypotheses:

\begin{prop} The varieties satisfying the hypotheses
of theorem 5.1 are $\emptyset\subset\pp 1$ and $\BA\pp 2\subset\ppp \JA$.
The varieties $X$ so produced are the varieties occuring as the
space of lines through a point of $G^{ad}$ where $G$ is an
exceptional group, i.e.,
$X$ is $v_3(\pp 1)$ and $G_w(\BA^3, \BA^6)$.\end{prop}

Note that the hypotheses force $I_2(Z)\simeq W$ as $\fh$-modules
because $Z$ is contained in the cubic hypersurface, whose equation
gives an equivariant inclusion $W\ra S^2W^*$ (by contraction).

\begin{theo} (algebraic version)
Let $\fl$ be reductive with one dimensional center and let  
 $W$ be an irreducible $\fl$-module, with a 
non trivial action of the center. Suppose that 
$\we^2W$ is irreducible, so that $\fg=\fl\op W\op W^*$ is a 
simple Lie algebra (see \cite{lm1}). Suppose moreover that $W$ is endowed
with
an  $\fl$-invariant cubic form, and that,
as an $\fl$-module,  $S^2W=W^*\op S$, with $S$ irreducible. 
Then $V=\CC\op W \op W^* \op\CC$ has a natural
structure of simple $\fg$-module. \end{theo}

We thus recover the   constructions of
Freudenthal  without  using division algebras. Moreover, our proofs
will show
that the constructions
work because of the irreducibility of $I_2( Z)$. This perspective
simplifies the computations. In the same spirit,  we construct below,
the invariant symplectic and quartic forms from
a unified perspective and without
use of composition algebras.

The equivalence of the two versions is as
follows: If $W\ne\CC$,  then $S= S\up 2W$
(the Cartan product of $W$ with itself). In general, if $Z\subset
\PP W$ is a closed orbit,  then $I_2(Z)$ is the complement to  
$S\up 2W^*$
in $S^2W^*$. \medskip

The theorem is proved in the same way as the results 
on minuscule varieties in \cite{lm1}
only the argument is simpler. The idea is to
define the natural action on each factor and to normalize
the actions such that the Jacobi identities hold.
\smallskip

To define the action of 
$\fg$ on $V$,   let  $C \in S^3W^*$ denote the   cubic and
$C^*\in S^3W$ denote the dual cubic.
Then $\fg = \fl\op W \op W^*$ acts on $V$ in the following way:
$\fh\subset\fl $ acts naturally on each factor, in particular trivially
on $\CC$ and $\CC^*$;
$\11\in\CC={\mathfrak z}(\fl)$ acts by multiplication by 
$-3/2, -1/2, 1/2, 3/2$ on the four respective components of $W$;
finally, the actions of $W$ and $W^*$ are given by the following 
formulae:
$$\begin{array}{rcccccccc}
t.(\a\op r\op s^*\op \b ^*) & 
= & 0&\op &3\a t&\op  & C(rt)& \op & \frac{1}{2}\langle t,s^*\rangle, \\
t^*.(\a\op r\op s^*\op \b ^*) & 
= & \frac{1}{2}\langle r,t^*\rangle &\op & C^*(s^*t^*) & \op &
3\b^*t^* & \op & 0. \end{array}$$

With this notation, the application
$\phi$ in \cite{lm1} in the special case of the 
minuscule theorem may be written as in \S 1.2.

The same construction in theorem 5.1   works to construct $Seg(\pp
1\times\pp 1\times
\pp 1) = G_w(\underline 0^3, \underline 0^6)$ out of
$\pp 0\sqcup\pp 0\sqcup\pp 0\subset\pp 2$, in fact, $Seg(\pp 1\times
Q^m)$ out of $\pp 0\sqcup Q^{m-2}$ where $Q^m$ is a  quadric hypersurface. 
The presence
of $\underline 0\pp 2= \pp 0\sqcup\pp 0\sqcup\pp 0\subset\pp 2$ should come
as
no suprise, as the Severi varieties also classify the
smooth connected base
schemes of the quadro-quadro Cremona transformations (see \cite{esb})
and $\pp 0\sqcup\pp 0\sqcup\pp 0\subset\pp 2$ is the base scheme
of the classical Cremona transform.
\smallskip

\begin{prop}
Let  $\Omega$ be the
  symplectic form   on 
$V=\BC\oplus W\oplus W^*\oplus\BC^*$   defined by
$$\O(\a\op r\op r^*\op \a ^*,\b\op s\op s^*\op \b^*)
=6(\a\b^*-\b\a^*)-(\langle r,s^*\rangle -\langle s,r^*\rangle ).$$

Then $\Omega$ is $\fg$-invariant.
\end{prop}

\proof The form $\O$ is clearly symplectic and $\fh$-invariant. Moreover,
if $u=\a\op r\op r^*\op \a^*$,   $v=\b\op s\op s^*\op \b^*\in V$,
 and $t\in W$, then
$$\begin{array}{rclcl}
\O (\11.u,v) & = & 
-9(\a\b^*+\b\a^*)+\frac{1}{2}(\langle r,s^*\rangle +\langle s,r^*\rangle ) 
& = & \O (\11.v,u), \\ 
\O (t.u,v) & = & 
-3(\langle t,s^*\rangle \a+\langle t,r^*\rangle \b)-C(rst) & = & \O
(t.v,u).
\end{array}$$
This means that $\O$ is $\CC$ and $W$-invariant, hence by symmetry
$W^*$-invariant as well. \qed


\subsection{The quartic invariant}

The five $\fg$-modules $V$ constructed above have a free invariant algebra,
generated in degree four (see e.g. \cite{brion}). 
We   write
down this quartic invariant in a unified way, in terms of the
$\fh$-invariant  cubic $C$ on $W$.
  
\begin{prop}\label{disc} The quartic polynomial defined 
for $w=\a\op r \op s^*\op\b^* \in V$ by 
$$Q(w)= (3\a\b^*-\frac{1}{2}\langle r,s^*\rangle )^2
+\frac{1}{3}(\b^*C(r^3)+\a C^*(s^{*3}))
-\frac{1}{6}\langle C^*(s^{*2}),C(r^2)\rangle $$
is a $\fg$-invariant form on $V$. 
\end{prop}

\proof $Q$ is obviously an $\fh$-invariant polynomial. It is also 
$\CC$-invariant: it is easy to check that each of its three terms 
$Q_1$, $Q_2$, $Q_3$ is $\CC$-invariant. Taking into account the 
symmetry of the expression of $Q$, we just need to check that it is 
invariant under the action of $W$, since it will immediately be
invariant also under the action of $W^*$. We compute the action of 
$t\in W$ on $Q_1$, $Q_2$, $Q_3$ separately:
$$\begin{array}{rcl}
t.Q_1(w) & = & -\frac{1}{2}
(3\a\b^*-\frac{1}{2}\langle r,s^*\rangle )C(r^2t), \\
t.Q_2(w) & = & \frac{1}{6}\langle t,s^*\rangle C(r^3)+
3\a\b^*C(r^2t)+\a C^*(s^{*2}C(rt)), \\
t.Q_3(w) & = & -\frac{1}{3}\langle C^*(s^*C(rt)),C(r^2)\rangle
-\a C^*(s^{*2}C(rt)).
\end{array}$$
It is then straightforward to check that the invariance of $Q$
is equivalent to the identity
$$2\langle C^*(s^*C(rt)),C(r^2)\rangle =\langle t,s^*\rangle C(r^3)+
3\langle r,s^*\rangle C(r^2t).$$

Let $\th :  W\ot W^*\ra \fg$ be the map dual to the action of $\fg$ on 
$W$. The fact that the Jacobi identities hold in $\fg$ amounts to the 
following lemma, which partly follows from \cite{lm1}, Proposition
5.1, and can be proved along the same lines. 

\begin{lemm}
We can normalize $\th, C$and $C^*$ in such a way that the 
following identities hold:
$$\begin{array}{rcl}
\th (r\ot t^*)s-\th (s\ot t^*)r & = & 
\langle s,t^*\rangle r-\langle r,t^*\rangle s, \\
\th (r\ot t^*)s+\th (s\ot t^*)r & = & 
2(\langle s,t^*\rangle r+\langle r,t^*\rangle s)-2C^*(t^*C(rs)).
\end{array}$$
\end{lemm}

Taking the difference of these two identities, and then
contracting with $C(r^2)$ gives precisely the equality we needed. \qed

\medskip 
The above expression  of the quartic invariant was rediscovered by
several authors in special cases \cite{igusa, fr}. 
One may wish to compare it with Freudenthal's uniform expression
for the quartic \cite{fr} p. 166 (which is only defined for
$V=\cZ_2(\BA)$).

\medskip 

For those who like to consider $D_4$ as an exceptional group, 
note that $\pp 1\times\pp 1\times\pp 1\subset\pp 7$, 
  arises as the tangent directions
to the lines through a
point of
$D_4^{ad}$.  In this case the quartic is the simplest instance
of Cayley's  
{\em hyperdeterminant}, see \cite{gkz}.

In particular, the above formula for the quartic applies
to induce the hyperdeterminant from the cubic on
$W= {\cal J}_3(\underline 0)= \BC^3$ defined by $C(a\oplus b\oplus c)=abc$.

   When we restrict 
the quartic form on $\cZ_2(\BA )$ to $\cZ_2(\underline 0  )$ we obtain the
hyperdeterminant and in turn, the hyperdeterminant determines a unique
$G$-invariant quartic form on each 
$\cZ_2(\BA )$. Moreover, specializing further we have:\smallskip

\begin{prop} The quartic   $Q$ on $\cZ_2(\BA )$ is the unique
 $Sp_6(\BA )$ invariant polynomial whose restriction to
the subalgebra $\Delta_{\cZ}\subset \cZ_2(\BA )$ 
is the classical discriminant.\end{prop}

Note that taking $\BA=\underline 0$, this gives a new characterization even
of the   hyperdeterminant.  \medskip

\begin{proof}
It is sufficient to show that the vector space
$\fsp_6(\BA)\Delta_{\cZ}$ is $\cZ_2(\BA )$.  
Suppose to the contrary that $\fsp_6(\BA)\Delta_{\cZ}=U$
is a proper subspace. Since each of the four \lq\lq matrix\rq\rq\ 
components of $\cZ_2(\BA )$ is weighted differently for the cubic
$C$, we see the subspace must be the sum of linear  subspaces
of each of the four components. In fact the two one-dimensional
components must be present as they are in $\Delta_{\cZ}$. Moreover,
the other two components are dual to one another so must be cut
equally. So it is sufficient to consider the action on
$\Delta_{\cal J}\subset\JA$. But the identity matrix is in an open orbit
and so we obtain everything.\end{proof}\qed

\medskip\noindent 
\begin{rema}The same construction works equally well with $Seg(\pp 1\times
Q^m)$,
to obtain a symplectic and quartic form
on $\BC^2\ot\BC^{m+2}$ from the cubic form
on $\BC\oplus \BC^m$, $C(a,b)=aq(b)$, where $\BC^m$ is equipped
with a quadratic form $q$.\end{rema}

\subsection{Orbits} 

Our description of the closed orbit $G_w(\BA^3, \BA^6)$ as the image of
$\phi$ 
has the following known consequence, which also follows from \cite{lm0} 
and \cite{lm1}, so   we omit the proof. 

Let $W$ be a vector space with 
a symplectic form $\o$. A  variety $Y\subset\ppp W$ is {\em
Legendrian} if for all
$y\in Y$, the  affine tangent space $\hat{T}_yY$ is a maximal 
$\o$-isotropic subspace.

\begin{prop} $G_w(\BA^3, \BA^6) \subset \ppp{\cal Z}_2(\BA)$
is a Legendrian variety.  Moreover, its tangential variety 
$\tau (G_w(\BA^3, \BA^6))$ is naturally isomorphic to its  dual variety
$G_w(\BA^3, \BA^6)^*
\subset \ppp{\cal Z}_2(\BA)^*$, and is the quartic hypersurface $(Q=0)$.
In other words, $\tau (G_w(\BA^3, \BA^6))=(G_w(\BA^3, \BA^6)_*)^*$.
\end{prop}

The isomorphism $\tau (G_w(\BA^3, \BA^6))=(G_w(\BA^3, \BA^6)_*)^*$
gives another connection between  the quartic invariant $Q$ and the theory
of
hyperdeterminants
\cite{gkz}. 

\medskip\noindent\begin{rema} It is unusual to have a smooth variety whose
dual has degree four and it would be interesting to classify such.
Zak \cite{zak} has classified the smooth varieties whose duals have degree
less than four; in degree two there is only the quadric hypersurface
and in degree three there are only ten examples,
$Seg(\pp 1\times \pp 2)$, its hyperplane section, the Severi varieties, and
the smooth  projections
of the Severi varieties. We are unaware of 
any   general method for constructing varieties with duals
of a given degree, but we record the following observation: \end{rema}

\begin{prop} The varieties of F-points in the
second row of the $n=d$ magic chart, namely
$v_2(\pp{d-1})$, $Seg(\pp{d-1}\times\pp{d-1})$ and $G(2, 2d)$
have dual varieties of degree $d$.\end{prop} 
 
\smallskip The orbit structure of each of the  varieties
$G_w(\BA^3,\BA^6)$ has 
already been studied  (e.g. in \cite{brion}), but their
  similarities seem  to have been overlooked. 
The following proposition follows from results in \cite{brion}.
We give two different short proofs along the lines of our study. 

\begin{prop} For each of the 
varieties of F-planes in the third row of the magic chart, there are
exactly four
orbits, the closures of which are ordered by inclusion:
$$
G_w(\BA^3, \BA^6)\subset\s_+(G_w(\BA^3, \BA^6))\subset\t (G_w(\BA^3,
\BA^6))\subset \PP V.$$ The equations of $\s_+(G_w(\BA^3, \BA^6))$
(respectively
$G_w(\BA^3, \BA^6)$) are given by  the first (respectively second)
derivatives
of the discriminant $Q$. The dimensions are respectively
$3m+3, 5m+3$ and $6m+6$.
\end{prop}

We also describe the intermediate orbit closure $\s_+(G_w(\BA^3, \BA^6))$:

\begin{prop}  $\s_+(G_w(\BA^3, \BA^6))$ can be described as 
\begin{enumerate}
\item the singular locus of $\t (G_w(\BA^3, \BA^6))$, 
\item the locus of points on a family of secant lines to $ G_w(\BA^3,
\BA^6)$
(a $(m+4)$-dimensional family for smooth points),
\item the locus of points on a  secant line to
$G_w(\BA^3,
\BA^6)$ isotropic for the symplectic form (unique if a smooth point),
in other words points on a secant line to two interwoven points
in the sense of Freudenthal, i.e. two points in a same
F-Schubert variety $\Sigma^{G_w(\BA^1,\BA^6), a}_{G_w(\BA^3,\BA^6) }$,
\item the locus of points on a tangent line to the distribution
$\tilde \sigma (\BA\pp 2)\subset TG_w(\BA^3,\BA^6)$.
\end{enumerate}
\end{prop}

\noindent\begin{rema} For the other Legendrian varieties that arise as
the space of lines through a point of an adjoint variety we
have the following orbit structures: 
for $v_3(\pp 1)$ (tangent directions to lines through a point of
$G_2^{ad}$), there are only three orbits as
$\s_+(v_3(\pp 1))=v_3(\pp 1)$. For
$\pp 1\times Q^m$  (tangent directions to lines through a point of
$SO(m)^{ad}$), the structure is the same as above except that
 $\s_+(\pp 1\times
Q^m)$ decomposes into two irreducible components, $Seg(\pp
1\times\pp{m+1})$
and $\{ [e\ot b + f\ot c] \mid
e,f\in\pp 1\ b\ww c\in G_o(2, \BC^{m+2})\}$, with the
exception of $Q^2=\pp 1\times\pp 1$ where there are three
components.
\end{rema}

\medskip

\noindent\begin{rema}
We determine the orbit structure using
an algorithm that is applicable in general.   The idea is to
infinitesimalize
the study and reduce the problem to a lower dimensional question. Let
$X=G/P\subset\ppp V$ be the closed orbit and fix $x\in X$. Then every $v\in
V$
is in some 
  $\hat T_x\up kX\backslash T_x\up{k-1}X$, where 
$\hat T_x\up kX$ denotes the $k$-th osculating 
space (see \cite{lm0}). Letting $H$ be a maximial semi-simple subgroup of
$P$,
each
$N_k =  \hat T_x\up kX/ T_x\up{k-1}X$ is an $H$-module and has
corresponding
 orbits. Say there are $p_k$ $H$-orbits in $N_k$
and $T\up d=V$.
Then there are at most $p_1+\hdots +p_d$ $G$-orbits and in fact there are
strictly less because different $H$-orbits will lead to the same
$G$-orbit.
\end{rema}

 \medskip\proofs 
 Write $G_w(\BA^3, \BA^6)=G/P$ and let
$Sl_3(\BA) \subset P$ be a maximal semi-simple subgroup. Then there are
four
$Sl_3(\BA)$-orbits in $  T_xG_w(\BA^3, \BA^6)$, namely $0$,  $\hat \BA\pp
2\backslash 0$,
$\hat
\s (\BA\pp 2)\backslash
\hat \BA\pp 2$, and $  T_xG_w(\BA^3, \BA^6)\backslash \hat \s(\BA\pp 2)$.
Since
$\BA\pp 2$ is the base-locus of the second fundamental form,  it gives the
same
$G$-orbit as $0$. Thus there are at most three $G$-orbits in $\t
(G_w(\BA^3,\BA^6))$. To see that there are indeed three, note that the
space $\s_+(G_w(\BA^3, \BA^6))$ of tangent directions to the 
distribution $\tilde\s (\BA\pp 2)$ 
is $G$-invariant, strictly  contains $G_w(\BA^3, \BA^6)$ (since the 
intersection of $G_w(\BA^3, \BA^6)$ with any of its tangent spaces is 
an $\hat \BA\pp 2$), and is properly contained in $\t
 (G_w(\BA^3,\BA^6))$ (its dimension being smaller). Using the 
rational map $\phi$ above, it is easy to check that the derivatives of
the quartic $Q$ vanish on $\hat\s (\BA\pp 2)\subset T_xG_w(\BA^3,
 \BA^6)$ for $x=(1,0,0,0)$, hence for any $x\in G_w(\BA^3, \BA^6)$. 
This implies that $\s_+(G_w(\BA^3, \BA^6))$ is the singular
locus of $\t (G_w(\BA^3, \BA^6))$. 

To prove that the equations of $G_w(\BA^3, \BA^6) $ are the second
derivatives 
of $Q$, we just notice that this space of quadratic equations define 
a non empty $G$-stable subset of $\PP V$, properly contained in 
the singular locus of $\t  (G_w(\BA^3, \BA^6))$. Because of the orbit
structure, this must be $ G_w(\BA^3, \BA^6)$. 
 Finally, since $\t (G_w(\BA^3, \BA^6))$ is a hypersurface, its 
complement  must be  an open orbit.

The second proposition follows
by observing that each of these characterizations defines a union of
orbits in
$\PP V$,  which is properly contained in $\t (G_w(\BA^3, \BA^6))$, 
but different from $G_w(\BA^3, \BA^6)$. Hence
each must coincide with $\s_+(G_w(\BA^3, \BA^6))$.\qed

\medskip

Our second proof  gives more information about the entry loci
and other geometric objects:

A generic point $p\in \t (G_w(\BA^3, \BA^6))$ lies on a unique
tangent line so it will be sufficient to show there are points
lying on a family of tangent lines but not on   $G_w(\BA^3, \BA^6)$.
Let $p$ be on a tangent line to $x\in G_w(\BA^3, \BA^6)$
such that $p$ corresponds to a vector $v\in T_xG_w(\BA^3, \BA^6)$
with the property that $[v]\in \t (\BA\pp 2)_x \subset
\ppp T_xG_w(\BA^3, \BA^6)$.
In this case there exists $y\in \BA\pp 2_x$ such that $v$ may also
be considered an element of $T_y\BA\pp 2_x$. Moreover,  we
may consider $y\in G_w(\BA^3, \BA^6)$ (as tangent directions
in $\BA\pp 2$ correspond to lines on $G_w(\BA^3, \BA^6)$).
On the other hand since $\t (\BA\pp 2)$ is degenerate,
there is an $\BA\pp 1=Q^m$'s of choice of $y$ (see \S 4 above), and
  any point   $z\in \BA\pp 1\subset G_w(\BA^3, \BA^6)$
  has a tangent vector $w\in T_zG_w(\BA^3, \BA^6)$ corresponding
to $p$. Thus $p$ is on an $(m+1)$-dimensional family of tangent
lines  
and thus an $(m+2)$-dimensional family of secant lines. By our
explicit description, the equivalence of 1,2 and 4   follows
and 3 follows from noticing that the $\BA\pp 2$ is an F-Schubert
variety associated to an F-point $a\in G_w(\BA^1,\BA^6)$.\qed

\medskip
Note that we recover  that there is no
  $G$-invariant polynomial on $V$  (up to constants) other than $Q$ and 
its powers.   
The geometric interpretation of $\s_+(G_w(\BA^3, \BA^6))$ has been
investigated 
in \cite{donagi} in the case   of the 
Grassmannian $G(3,6)\subset\ppp (\La 3 \BC^6)$.

 Since there are only four orbits, we also have:

\begin{prop} With the notations   above,   $\s_+(G_w(\BA^3, \BA^6))$
is self-dual. 
\end{prop}
 
 The following
proposition can be proved in the same way that Proposition 3.2. 
 
\begin{prop}  $\s_+(G_w{(\BA^3, \BA^6)})$ is ruled by 
the  $\pp{m+3}$'s that are the linear spans of the F-Schubert
varieties 
$$\Sigma_{G_w(\BA^3, \BA^6)}^{Sp_6(\BA)^{ad}}\simeq Q^{m+2}.
$$  In particular, 
a smooth point of $\s_+(G_w{(\BA^3, \BA^6)})$  lies on a unique
$\pp{m+3}=< Q^{m+2} >$.
\end{prop}

\noindent\begin{rema} 
Let $V$ be an irreducible $\fg$-module.
The decomposition of $\fg$-modules 
$$
S^2V=S\up 2 V\oplus W\oplus stuff
$$
with $W$ irreducible implies that the closed $G$-orbit
in $\ppp W$ induces a variety of quadrics of constant rank
on $V$, and linear spaces on the closed orbit furnish
linear systems of quadrics.  

Linear systems of quadrics of constant rank arise as the second fundamental
forms of degenerate dual varieties
(see \cite{IL}), and few examples of such systems
(or smooth varieties with degenerate duals) are known.

The interpretation of F-Schubert varieties associated to
planes  as a family of
quadrics  on $G_w(\BA^3,\BA^6)$  is related to the  
decomposition of the symmetric square of $\cZ_2(\BA) $ as 
$$S^2\cZ_2(\BA)=S^{(2)}\cZ_2(\BA)\op\fs\fp_6(\BA).$$
An  element $X\in\fs\fp_6(\BA)$ defines a quadratic form on $\cZ_2(\BA)$,
namely $q_X(u)=\om (Xu,u)$. 
  In particular, we obtain   varieties
of quadrics of constant rank. \end{rema}

\begin{prop} Let $m=1,2,4,8$. The adjoint variety $Sp_6(\BA)^{ad}$,
parametrizes  a variety of dimension $4m+1$ of quadrics of 
rank $m+4$ on $\cZ_2(\BA)=\BC^{6m+8}$.
\end{prop}

One can take linear spaces on these varieties on these varieties
(except for $v_2(\pp 5)$) to get linear systems of quadrics of constant 
rank.

\section{Adjoint varieties of the exceptional groups}

As shown in \cite{lm1}, $\fe (\BA)$ may be constructed
from $\fsp_6(\BA)$, and the adjoint variety $X^{ad}_{E(\BA )}$
may be constructed from $G_w(\BA^3, \BA^6)$ via 
a rational map of degree four, given in terms of the 
quartic invariant $Q$ on $V$.
The construction  also reproduces the five step $\ZZ$-grading
of $\fe (\BA)$ and the filtration of $\fe (\BA)$ induced
by $U(\fs\fp_6 (\BA))$. 
Continuing the notation of the previous section,
$$
\fe (\BA ) = \BC^*\oplus V^*\oplus (\BC\oplus \fsp_6(\BA))
\oplus V
\oplus \BC
$$
where using the composition algebra model, $V={\cal Z}_2(\BA)$.
The adjoint variety $E(\BA)^{ad}$ is the image of the rational
mapping $\psi$ described in \S 1.2.

\medskip

While $\ppp (\fe (\BA))$ does not have a finite number of
$E(\BA )$-orbits, there are only a finite
number of $E(\BA )$-orbits in $\s ( E(\BA )^{ad})\subset \ppp (\fe (\BA))$
which we now describe.
  The following theorem improves upon recent 
results in \cite{kaji} where it is shown that $\s (G^{ad})$ contains
an open orbit for any simple group $G$.
We show that it is actually the union of a finite number of orbits
and exhibit them explicitly:

\begin{theo} Let $\fe (\BA)$ respectively
denote   $\ff_4, \fe_6,
\fe_7,\fe_8$, and  $m=1,2,4,8$. Let $E(\BA)^{ad} \subset\ppp (\fe (\BA ))$
denote the adjoint variety, the closed
$E (\BA)$-orbit. Then $\s (E(\BA)^{ad})=\tau (E(\BA)^{ad})$
  and the $E(\BA)$-orbit closures
in $\s (E(\BA)^{ad})$ are as follows:
$$E(\BA)^{ad}\subset\s_{(2m+7)}(E(\BA)^{ad})  \subset\s_{(3)}(E(\BA)^{ad}) 
\subset\s_{(1)}(E(\BA)^{ad}) 
\subset\s(  E(\BA)^{ad}  ).$$

These orbits are respectively of dimensions 
$6m+9, 10m+11, 12m+15,
 12m+17$ and $12m+18$.

 Moreover, the open orbit in $\s (E(\BA)^{ad})$ is   a 
semi-simple orbit, while the four others are projectivizations 
of nilpotent orbits. 
\end{theo}

With this notation, the orbit closure $\s_{(m)}(E(\BA)^{ad})$
has codimension $m$ inside $\s (E(\BA)^{ad})$,
so a general point of $\s_{(m)}(E(\BA)^{ad})$ 
has an $(m+1)$-dimensional entry locus.

\begin{prop} 
The  orbit closures above can be described as follows:
\begin{itemize}
\item $\s_{(1)}(E(\BA)^{ad})=\s (E(\BA)^{ad})\cap Q_{Killing}$, 
 where $Q_{Killing}$ is the quadric hypersurface defined by the
 Killing form. Equivalently,  it is   the closure of the orbit of
points belonging to  a  unique  tangent line to the distribution 
$T_1E(\BA)^{ad}$ of contact hyperplanes in $TE(\BA)^{ad}$, and the
points on a secant line of two {\it hinged} points in
the sense of Freudenthal (see below).
\item $\s_{(3)}(E(\BA)^{ad})$ is the closure of the orbit consisting 
of points on  a tangent line  to the distribution 
$\tilde \tau (G_w(\BA^3,\BA^6))\subset TE(\BA)^{ad}$, a one 
dimensional family of such. 
\item $ \s_{(2m+7)}(E(\BA)^{ad})$ is the closure of the orbit
  consisting of points belonging to an $(m+2)$-dimensional family of 
tangent lines to the distribution $\tilde \s_+ (G_w(\BA^3,\BA^6))
\subset TE(\BA)^{ad}$,  equivalently of points on a secant line of two   
{\it interwoven}  points in the sense of Freudenthal. 
\end{itemize}
\end{prop} 
 
Note that $ \s_{(3)}(E(\BA)^{ad}) $ cannot be detected
from Freudenthal's geometries. Thus, as in \cite{lm0}, the
perspective
of Freudenthal  and Tits is extremely useful for understanding the
projective geometry, but it does not reveal the full story.
\medskip

\begin{coro} $\s_{(2m+7)}(E(\BA)^{ad}) $ is ruled by 
the  $\pp{m+5}$'s that are the linear spans of the F-Schubert
varieties 
$$\Sigma^{X_{F-points}^{E(\BA)} }_{E(\BA)^{ad}}\simeq Q^{m+4}.
$$
In particular, 
a smooth point of $\s_{(2m+7)}(E(\BA)^{ad})$  lies on a unique
$\pp{m+5}=< Q^{m+4}>$.
\end{coro}

\medskip
Before entering into the proof of the theorem, we recall
Freudenthal's incidence relations for points of F-symplecta,
that is points of $E(\BA )^{ad}$. They can be:

 {\em joined}, which means they are contained in
a unique F-plane. In other words, two points $x,y\in E(\BA )^{ad}$ are
joined if
their  secant line $\pp 1_{xy}$ is contained in $E(\BA )^{ad}$;
  
{\em interwoven}, which means they intersect in an F-point.
In other words, two points $x,y\in E(\BA )^{ad}$
are  interwoven if they are contained in a 
$\Sigma^{X_{F-points}}_{E(\BA)^{ad}}= Q ^{m+4}
\subset E(\BA )^{ad}$; if this F-Schubert variety is not 
unique, then $x$ and $y$ are joined;
 
{\em hinged}, which means they are joined to a third F-symplecton. In other
words, 
$x,y\in E(\BA )^{ad}$ are (strictly) hinged if there exists
(a unique) $z\in E(\BA )^{ad}$
such that the secant lines $\pp 1_{zx}, \pp 1_{zy}$ are contained in
$E(\BA )^{ad}$.
Equivalently,  $x,y\in E(\BA )^{ad}$ are   hinged if
their secant line is
contained in the quadric hypersurface defined by the Killing 
form;
 
{\em generic}, i.e.,  not hinged.  
\medskip

\proofs By   \cite{kaji},   
$\s (E(\BA)^{ad})= \t (E(\BA)^{ad})$ and has secant defect one.
Also, as noted in \cite{kaji}, fixing a Cartan subalgebra and a set
of simple roots for $\fe (\BA)$, the orbit of $X_{\tilde{\a}}
+X_{-\tilde{\a}}$ is open in $\s (E(\BA)^{ad})$, where $\tilde{\a}$ 
denotes the maximal root. 
This element is   semi-simple, and conjugate to a 
multiple of $H_{\tilde{\a}}$. This proves that the open orbit 
is isomorphic to the (semi-simple) orbit of $H_{\tilde{\a}}$.
Consider in $\fe (\BA)$ the cone $C$ over this orbit, and an element $x$ of
its closure. The semi-simple part of $x$ must be conjugate to
$\l H_{\tilde{\a}}$ for some scalar $\l$. If $\l\ne 0$ and $x$ is not 
semi-simple, the cone over the orbit it generates is of dimension 
strictly bigger than $C$, which is absurd. Hence $x$ is semi-simple 
or nilpotent, and this proves that $\overline{C}/C$ is  a union of 
nilpotent orbits.  

 There potentially are
four kinds of elements of $\t (E(\BA)^{ad})$, corresponding to
$v\in T_xE(\BA)^{ad}\backslash T_1$,  
$v\in T_{1x}\backslash \s G_w(\BA^3,\BA^6)_x$,
$v\in \s (G_w(\BA^3,\BA^6))_x\backslash   G_w(\BA^3,\BA^6)_x$
and $v\in   G_w(\BA^3,\BA^6)_x$. The last type is the same as
a point on $E(\BA)^{ad}$.
Let $p\in \t (E (\BA)^{ad})$.

 Consider the case where there exists an $x\in E (\BA)^{ad}$ and
 $v\in T_{1x}E (\BA)^{ad}$ with $p$ on the line corresponding to $v$.
Since $\s (G_w(\BA^3,\BA^6))_x= \ppp T_{1x}$, there exist
$y,z\in  G_w(\BA^3,\BA^6)_x\subset E (\BA)^{ad}$ such that
$p\in \pp 1_{yz}$. Since   $\pp 1_{xz},\pp 1_{xy}\subset
E (\BA)^{ad}$, we see that $p$ is indeed on a secant line
of two hinged points, showing the equivalence of the first and
third characterizations, modulo the unicity in the
second: but this follows from Freudenthal's remark that if 
two F-symplecta are multiply hinged, that is joined to several others 
F-symplecta, they must be interwoven. Finally, the third and second 
characterizations are equivalent from Freudenthal's observations
again. 

\smallskip 
Now consider the case where there exists an $x\in E (\BA)^{ad}$ and
 $v\in \hat \t (G_w(\BA^3,\BA^6)_x)\subset T_{1x}E (\BA)^{ad}$ with $p$ on
the
line corresponding to
$v$. Let $y\in  G_w(\BA^3,\BA^6)_x$ be the (in general unique)
point such that $v$ corresponds to a vector in $T_yG_w(\BA^3,\BA^6)_x$.
By the same argument as in \S 5.3 above, $p$ lies on a tangent line
to all $z\in \pp 1_{xy}$ and moreover these tangent lines
are tangent to the distribution $\tilde \t G_w(\BA^3,\BA^6)$.
Moreover, for each $x$, $y$ and $v$ are uniquely determined, so the
component
of the   locus tangent to $\tilde \t (G_w(\BA^3,\BA^6))$ passing through
$x$ is a $\pp 1$. The dimension count follows.

\smallskip 
Finally consider the case where there exists an $x\in E (\BA)^{ad}$ and
 $v\in \hat \s_+ G_w(\BA^3,\BA^6)_x\subset T_{1x}E (\BA)^{ad}$ with $p$ on
the
line corresponding to
$v$. Now
there exist $y,z\in   \hat G_w(\BA^3,\BA^6)_x \subset E(\BA)^{ad}$
such that $p\in \pp 1_{yz}$, in fact a 
$\hat Q^{m+3}\subset 
\hat G_w(\BA^3,\BA^6)_x 
$   of such points so this orbit closure does not coincide with
any of the others.  The points $y,z$ are interwoven as they are
both contained in a $Q^{m+4}\subset E(\BA)^{ad}$,  that
is an F-Schubert variety  
$$\Sigma^{X^{\fe (\BA)}_{F-points} }_{E(\BA)^{ad}}$$
 showing the equivalence of the two characterizations.
Each of these F-Schubert varieties generates a $\pp{m+5}$ in 
$\PP(\fe (\BA))$. 

We prove that two such $\pp{m+5}$'s,  if they are not equal, 
can intersect only inside ${E(\BA)^{ad}}$. Suppose the contrary, and
take a generic line $l$ in their intersection. It cuts ${E(\BA)^{ad}}$
exactly in two points $u$ and $v$, since the intersection of
${E(\BA)^{ad}}$
with each of our $\pp{m+5}$'s is a quadric. But then $u$ and
$v$ are doubly interwoven, hence joined, thus the line $l$ is contained
in ${E(\BA)^{ad}}$, a contradiction. 

This proves that a generic point of $\s_{(2m+7)}(E(\BA)^{ad})$ belongs
to a unique $\pp{m+5}$ generated by an F-Schubert variety. The dimension
follows because $\tdim X^{\fe (\BA)}_{F-points} + (m+5) = 10m+11$.
\qed

\medskip One can check from the tables in \cite{col} that there exists 
nilpotent orbits with the dimensions claimed in the proposition (minus
one, because of the projectivization). If we exclude $F_4$, they are
respectively labelled $A_1$, $2A_1$, $3A_1$ ($3A_1'$ in the case of 
$E_7$), and $A_2$.

\medskip\rem The symmetric squares of the exceptional 
simple Lie algebras $\fe (\BA)$  have a uniform 
decomposition into irreducible components: 
$$S^2\fe (\BA)=S^{(2)}\fe (\BA) \op W \op \CC,$$
where $S^{(2)}\fe (\BA)$ denotes the Cartan product of $\fe (\BA)$ with
itself, 
the $\CC$ component is given by the Killing form, and the other 
component $W$ is the ambient space for F-points. 
In particular, we obtain:

\begin{prop} $X^{\fe (\BA)}_{F-points}$ parametrizes a
variety of dimension $9m+6$ of quadrics of rank $m+6$ on $\fe (\BA)$.
\end{prop}

We now describe the orbit structure of $\s (G^{ad})$
for the remaining simple groups.
In each case the following properties hold:  the open orbit
is semi-simple and the others are nilpotent;
 there is the orbit
$\s_{(1)}(G^{ad})=
\s (G^{ad})
\cap Q_{Killing}$, equivalently, the points on a tangent line
of the distribution of contact hyperplanes $T_1G^{ad}$;
and $\tdim \s
(G^{ad})= 2\tdim G^{ad}$,  so  $\tdim \s_{(m)}(G^{ad})= 2\tdim G^{ad}-m$.

Note that the orbit structure  
   is not as uniform for the  classical 
groups as for the exceptional groups. 

\medskip\noindent ${\bf G= G_2}$. 
The orbit structure inside the secant variety is 
$$G_2^{ad} \subset \s_{(3)}(G_2^{ad})
\subset \s_{(1)}(G_2^{ad})\subset \s (G_2^{ad})= \t(G_2^{ad}).
$$
$\tdim G_2^{ad}= 5$. The
open orbit is semi-simple and the others are nilpotent.
In particular, $\s_{(1)}(G_2^{ad})$ is the closure of the 
projectivization of the subregular nilpotent orbit. 
 The orbit closure
$\s_{(3)}(G_2^{ad})$ consists of points on a tangent line
to the distribution $\tilde \t (v_3(\pp 1))$. 

\medskip\noindent ${\bf G=SL_n}$. 
Note that $SL_n^{ad}=\FF_{1,n-1}\subset\ppp
(\fsl_n)$ is  a partial flag variety  of  dimension $2n-3$. 
Here the orbit structure inside the secant variety is 
$$SL_n^{ad}=\FF_{1,n-1}\subset \s_{(3)}(\FF_{1,n-1})
\subset \s_{(1)}(\FF_{1,n-1})\subset \s(\FF_{1,n-1})= \t(\FF_{1,n-1}).
$$  
The intermediate orbits correspond to endomorphisms that in Jordan normal
form consist of one $3\times 3$ nilpotent block, 
and two $2\times 2$ nilpotent blocks respectively.  
 
\medskip\noindent ${\bf G=SO_n}$. 
Here $SO_n^{ad}=G_o(2,n)\subset\ppp (\La 2\CC^n)$ is an isotropic
Grassmanian
for a quadratic form $Q$. The orbit structure is as follows:
$$SO_n^{ad}=G_o(2,n)\subset \matrix \s_{(2n-9)}(G_o(2,n))
\\ \s_{(7)}(G_o(2,n))\endmatrix 
\subset \s_{(1)}(G_o(2,n))\subset \s (G_o(2,n)).
$$
The orbit closures are {\sl not} totally ordered by
inclusion. This orbit structure is not suprising because
$\tbase II= \pp 1\times \tilde Q$ and
$\s_+(\pp 1\times \tilde Q)$ has two irreducible components.
The orbit closures  $\s_{(2n-9)}(G_o(2,n))$
and $\s_{(7)}(G_o(2,n))$ correspond to points on a
tangent line to one of the two corresponding distributions.
(Here, if $E\in G_o(2,n)$, then $\tilde Q\subset\ppp E\upperp/E$
is a quadric hypersurface.)
Note that, unlike in other cases,  if a point lies on a tangent line to the
distribution $\tilde\t (\pp 1\times \tilde Q)$, it is automatically
also on a tangent line to the distribution $\tilde\s_+
(\pp 1\times \tilde Q)$.

We may see the orbits from the global
geometry as follows:  let $P,P'\in G_o(2,n)$ be distinct points.  
Let 
$M= P+P'\subset \CC^n$. Then $\tdim M=3$ or $4$. If $\tdim M=3$,
then $\trank Q|_M= 0$ or $1$.
 When $\trank Q|_M= 0$,  $P$ and $P'$ are perpendicular
and the corresponding secant is contained in $G_o(2,n)$.
When $\trank Q|_M= 1$,   the corresponding orbit is $\s_{(2n-9)}
(G_o(2,n))$.  If 
$\tdim M=4$,  $\trank Q|_M= 0$, $2$ or $4$. These cases determine
orbit closures as follows: Rank zero occurs when $P$ and
$P'$ are perpendicular; the corresponding orbit 
is an open subset of a $\PP^5$-bundle over $G_o(4,n)$, its closure 
is $\s_{(7)}(G_o(2,n))$. Rank two occurs when $P'$ contains a line
perpendicular to $P$, the corresponding orbit 
is an open subset of a $G(2,n-4)$-bundle over $G_o(2,n)$; 
its closure is $\s_{(1)}(G_o(2,n))$. Rank four is the generic case.

\medskip\noindent ${\bf G=Sp_{2n}}$. 
Here $Sp_{2n}^{ad}=v_2(\PP^{2n-1})\subset\ppp S^2\CC^{2n}$. 
Taking two distinct lines $l$ and $l'$, the plane they generate is
either isotropic or not. The orbit corresponding to the isotropic case 
is an open subset of a $\PP^2$-bundle over $G_w(2,2n)$. 
The orbit structure of the secant variety of $Sp_{2n}^{ad}\subset
\ppp (\fsp_{2n})$  is therefore
$$Sp_{2n}^{ad}=v_2(\PP^{2n-1})
\subset \s_{(1)}(Sp_{2n}^{ad})\subset \s (Sp_{2n}^{ad}).
$$
  This is the most degenerate case.

\section{Desingularizations}

Let  $X\subset\BP V$ be a smooth variety. If  
the tangential variety of $X$, $\tau (X)\subset\BP V$ is nondegenerate,
then it admits a desingularization 
$ \tilde TX \ra \tau (X) $
where $\tilde TX$ is the bundle of embedded tangent projective spaces.
Similarly, if the dual variety $X^*\subset\ppp V^*$ is nondegenerate,
it admits a desingularization $\ppp N^*\ra X^*$
where $N^*$ denotes the conormal bundle of $X$.

When $X$ is homogeneous, both of these desingularizations
are examples of what Kempf \cite{kempf} calls the {\em collapsing of a
vector bundle}.   In particular, whenever $\t (X)$ or $X^*$ is
nondegenerate,
it has rational singularities and one can explicitly describe
its desingularization via Tits transforms.

\subsection{Orbits in $\PP(\JA)$}

 Here since $\s (\BA\pp 2)\simeq (\BA\pp 2_*)^*$   the above
discussion   applies to desingularize $\t (\BA\pp 2)= \s (\BA\pp 2)$.
Moreover, the  bundle $M =N^* (\AA\PP^2_*)(-1)$ can be 
described as follows: its fiber over an F-line is the linear 
subspace of $\JA$ generated by the F-Schubert variety consisting 
of F-points incident to this F-line. 

\begin{prop} Let $G/P=\BA\pp 2\subset\ppp V$ be a Severi variety. Let 
$\BA\pp 2_*\subset\ppp V^*$ denote the Severi variety in the dual
projective space, and $M$ be as above. There is a natural diagram
$$\begin{array}{ccc}
 E=\tilde Q^m & \longrightarrow & \BA\pp 2 \\ \cap & & \cap \\
 \PP M & \stackrel{f}{\longrightarrow} & \s (\BA\pp 2) \\
 \pi\downarrow & & \\ \BA\pp 2_* & & 
\end{array}$$
where $f$ is a desingularization of $\s (\BA\pp 2)$. The exceptional
divisor $E$ of $f$ is naturally identified with the $G$-homogeneous 
space consisting of pairs of incident F-points and F-lines, 
with its two natural projections over $\BA\pp 2$ and $\BA\pp 2_*$.
In other words, $E$ is the set of points in $(\BA\pp 2_*)^*$
tangent to $ \BA\pp 2_*$  along a quadric $\AA\PP^1\simeq Q^m$.   
\end{prop}

\noindent In the four diagrams below, 
we indicate  the nodes defining
the space of F-points $\BA\pp 2$
with   black dots  and those defining the F-lines
$\BA\pp 2_*\simeq\BA\pp 2$  with stars. 
The bundle $M$ on $\BA\pp 2_*$ (which is defined by 
the same node as $\AA\PP^2$),  
and the quadric inside the fibers of $\PP M$ are below the diagrams. 

\medskip
\begin{center}\begin{tabular}{cccc} 
$v_2(\PP^2)$ & $\PP^2\times \PP^2$ & $G(2,6)$ & $\hspace*{1mm}\OO\PP^2$ \\

\setlength{\unitlength}{2.5mm}
\begin{picture}(10,3)(-2.5,-1)
\put(3,0){$\circ$}
\put(1,0){$\bullet$}
\put(3,0){$*$}
\put(1.55,.35){\line(1,0){1.55}}
\end{picture} & 

\setlength{\unitlength}{2.5mm}
\begin{picture}(10,3)(0,-1)
\multiput(3,0)(2,0){2}{$\circ$}
\multiput(1.55,.35)(4,0){2}{\line(1,0){1.55}}
\put(1,0){$\bullet$}
\put(7,0){$\bullet$}
\put(3,0){$*$}
\put(5,0){$*$}
\end{picture} &

\setlength{\unitlength}{2.5mm}
\begin{picture}(10,3)(0,-1)
\multiput(0,0)(2,0){5}{$\circ$}
\multiput(0.55,.35)(2,0){4}{\line(1,0){1.55}}
\put(2,0){$\bullet$}
\put(6,0){$*$}
\end{picture} & 
 
\setlength{\unitlength}{2.5mm}
\begin{picture}(11,3)(0,-1)
\multiput(0,0)(2,0){5}{$\circ$}
\multiput(0.5,.3)(2,0){4}{\line(1,0){1.6}} 
\put(0,0){$\bullet$}\put(8,0){$*$}
\put(4,-2){$\circ$}
\put(4.3,-1.4){\line(0,1){1.5}}\end{picture} \\ 

$S^2\CC^2$ & $\CC^2\ot\CC^{2*}$ & $\Lambda^ 2Q$ & $\hspace*{5mm}\CC^{10}$
\\
$v_2(\PP^1)$ & $\PP^1\times\PP^1$ & $G(2,4)$ & $\hspace*{5mm}\QQ^8$
 \end{tabular}\end{center}\medskip

\subsection{Orbits in $\PP ({\cal Z}_2(\BA ))$}

The remarks at the beginning of
this section apply to $\t (G_w(\BA^3, \BA^6))$, which is singular 
exactly along $\s_+(G_w(\BA^3, \BA^6))$. 
It can be desingularized by a collapsing which
 answers a question of Kempf in the case 
of $E_7$, who   failed to
observe the orbit corresponds to a nondegenerate tangential variety. 

\begin{prop} Let $\tilde{T}G_w(\BA^3, \BA^6)$ be the  
bundle of embedded tangent spaces of   $G_w(\BA^3, \BA^6)$, 
whose associated vector bundle has rank $\tdim G_w(\BA^3, \BA^6)+1$.  There
is a
natural diagram
$$\begin{array}{ccc}
 E=\tilde\sigma(\BA\pp 2) & \longrightarrow & \s_+(G_w(\BA^3, \BA^6)) \\
\cap & &
\cap
\\
  \tilde{T}G_w(\BA^3, \BA^6)  & \stackrel{g}{\longrightarrow} & \t
(G_w(\BA^3, \BA^6)) \\
 \pi\downarrow & & \\ G_w(\BA^3, \BA^6)
 & & 
\end{array}$$
where $g$ is a desingularization of $\t (G_w(\BA^3, \BA^6)
)$. 
\end{prop}

The exceptional divisor 
$E$ is singular and it too can be desingularized by
a homogeneous projective bundle over the $G$-homogeneous space of 
pairs of incident F-points and F-planes. 
The singular locus of $E$ is
$\tilde{\BA\pp 2}\subset TG_w(\BA^3\BA^6)$ and its fibers over 
$G_w(\BA^3\BA^6)$ are cones over the Severi varieties $\AA\PP^2$.
Outside this locus, $E$ is a $Q^{m+1}$-bundle: a generic point $p$ in 
$\s_+(G_w(\BA^3, \BA^6) )$ is contained in the linear span $\PP^{m+3}$ 
of a unique F-Schubert variety
$\Sigma^{Sp_6(\BA)^{ad}, y}_{G_w(\BA^3,\BA^6)}$, which is a quadric
$Q^{m+2}$ 
inside this $\PP^{m+3}$. The fiber $g^{-1}(p)$ is then the section of 
$Q^{m+2}$ by the hyperplane perpendicular to $p$. 

The orbit closure   $\s_+(G_w(\BA^3, \BA^6) )$ also admits a natural
desingularization by a collapse given by
Freudenthal geometry. 
Let  $S$ be the homogeneous vector bundle on 
$Sp_6(\BA)^{ad}$ (the space  of F-points)   defined by the node of the
Dynkin diagram corresponding to F-planes, i.e., 
the bundle whose fiber at $y$ is the linear span of the  F-Schubert
variety
$\Sigma^{Sp_6(\BA)^{ad}, y}_{G_w(\BA^3,\BA^6)}$. 

\begin{theo} There is a natural diagram
$$\begin{array}{ccc}
 E=\tilde Q^{m+2} & \longrightarrow & G_w(\BA^3, \BA^6) \\ \cap & & \cap \\
 \PP S & \stackrel{f}{\longrightarrow} & \s_+(G_w(\BA^3, \BA^6)) \\
 \pi\downarrow & & \\ Sp_6(\BA)^{ad}= G_w(\BA^1,\BA^6) & & 
\end{array}$$
where $f$ is a desingularization of $\s_+(G_w(\BA^3, \BA^6))$. The
exceptional
divisor $E$ of $f$ is naturally identified with the $G$-homogeneous 
space consisting of pairs of incident F-points and F-planes, 
with its two natural projections over $G^{ad}$ and $G_w(\BA^3, \BA^6)$.  
The intersection of this divisor with a fiber of $\pi$ is a 
quadratic hypersurface. 
\end{theo}

 \noindent Our four examples of the above situation are the 
following, where we indicate   the nodes defining
the space of F-points with black dots, and  those
defining the F-planes with stars.  
The vector bundle $S$  and the quadric inside the fibers of $\PP S$
are given below the diagrams. 

\medskip
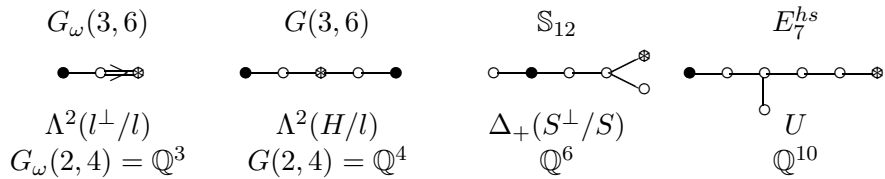
\begin{center}\begin{tabular}{cccc} 
$G_{\o}(3,6)$ & $G(3,6)$ & $\SS_{12}$ & $E_7^{hs}$ \\
\setlength{\unitlength}{2.5mm}
\begin{picture}(10,3)(-2.5,-1)
\multiput(0,0)(2,0){3}{$\circ$}
\put(4,0){$*$}
\put(0.55,.35){\line(1,0){1.55}}
\multiput(2.55,.25)(0,.2){2}{\line(1,0){1.55}} 
\put(2.75,0){$>$}
\put(0,0){$\bullet$}\end{picture} & 
\setlength{\unitlength}{2.5mm}
\begin{picture}(10,3)(0,-1)
\multiput(0,0)(2,0){5}{$\circ$}
\multiput(0.55,.35)(2,0){4}{\line(1,0){1.55}}
\put(0,0){$\bullet$}
\put(8,0){$\bullet$}
\put(4,0){$*$}
\end{picture} & 
\setlength{\unitlength}{2.5mm}
\begin{picture}(10,3)(1,-1)
\multiput(2,0)(2,0){4}{$\circ$}
\multiput(2.55,.35)(2,0){3}{\line(1,0){1.55}} 
\put(8.5,.35){\line(2,-1){1.55}}
\put(8.5,.35){\line(2,1){1.55}}
\put(10,.9){$\circ$}
\put(10,-.9){$\circ$}
\put(10,.9){$*$}
\put(4,0){$\bullet$}\end{picture} &
\setlength{\unitlength}{2.5mm}
\begin{picture}(11,3)(0.8,-1)
\multiput(0,0)(2,0){6}{$\circ$}
\multiput(0.5,.3)(2,0){5}{\line(1,0){1.6}} 
\put(0,0){$\bullet$}\put(10,0){$*$}
\put(4,-2){$\circ$}
\put(4.3,-1.4){\line(0,1){1.5}}\end{picture} \\ 
$\Lambda^ 2 (l^{\perp}/l)$ & $\Lambda^ 2 (H/l)$ & 
$\Delta_+(S^{\perp}/S)$  & $U$ \\
$G_{\o}(2,4)=\QQ^3$ &  $G(2,4)=\QQ^4$ & $\QQ^6$ & $\QQ^{10}$
 \end{tabular}\end{center}\medskip

\subsection{Orbits in  $ \ppp (\fe (\BA) ) $}

First note that $\tilde TE(\BA)^{ad} \ra \s (E(\BA)^{ad})$
provides a desingularization of $\s (E(\BA)^{ad})$  as $\s (E(\BA)^{ad})$
coincides with the tangential variety which is nondegenerate. The
desingularization is as follows:

\begin{prop} Let $E(\BA)^{ad}$ be a variety of F-symplecta
and let  $  {T_1}E(\BA)^{ad}\subset   TE(\BA)^{ad}$  be the 
bundle of contact hyperplanes.     
   There is a natural diagram
$$\begin{array}{ccc}
 \tilde \t G_w(\BA^3,\BA^6)  & \longrightarrow & \s_{(3)}(E(\BA)^{ad})\\
\cap & &\cap \\
  \tilde{T_1}E(\BA)^{ad} & \longrightarrow
&\s_{(1)}(E(\BA)^{ad})\\
\cap & & \cap \\
  \tilde{T }E(\BA)^{ad} & \stackrel{g}{\longrightarrow} &\s (E(\BA)^{ad})\\
 \pi\downarrow & & \\ E(\BA)^{ad}
 & & 
\end{array}$$
where $g$ is a desingularization of $\s(E(\BA)^{ad})$. 
\end{prop}

 Note that since $\s(E(\BA)^{ad})$ is normal (being the image of a
 collapsing), it is smooth in codimension one, hence the open orbit 
in the hypersurface 
$\s_{(1)}(E(\BA)^{ad})$ is contained in the smooth locus
of $\s(E(\BA)^{ad})$. This
implies that $g$ above is also a desingularization of
$\s_{(1)}(E(\BA)^{ad})$. In particular, we recover the fact that 
a general point of 
$\s_{(1)}(E(\BA)^{ad})$ belongs to a unique tangent line to the 
distribution $\tilde{T_1}E(\BA)^{ad}$.

Finally, here is a desingularization of $\s_{(2m+7)}(E(\BA)^{ad})$. 
Let $S$ be the bundle on $X^{E(\BA)}_{F-points}$ induced by the F-Schubert
varieties $\Sigma^{X^{E(\BA)}_{F-points}}_{E(\BA)^{ad}}\simeq Q^{m+4}$.

\begin{prop} Let $E(\BA)^{ad}$ be a variety of F-symplecta.     
   There is a natural diagram
$$\begin{array}{ccc}
 E =\tilde Q^{m+4}  & \longrightarrow &  E(\BA)^{ad}\\
\cap & &
\cap
\\
  \ppp S & \stackrel{h}{\longrightarrow} &\s_{(2m+7) }(E(\BA)^{ad})\\
 \pi\downarrow & & \\ X^{E(\BA)}_{F-points}
 & & 
\end{array}$$
where $h$ is a desingularization. 
\end{prop}

 \noindent Here the four examples of the above situation are the 
following, where we indicate   the nodes defining
the space of F-points with black dots, and  those
defining the vector bundle $S$ on $X^{E(\BA)}_{F-points}$
with stars.  

\medskip
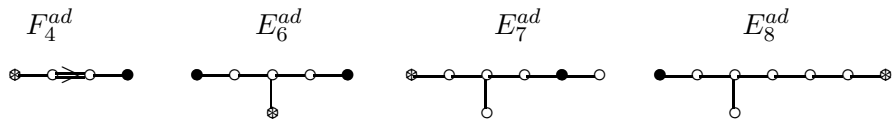
\begin{center}\begin{tabular}{cccc} 
$F_4^{ad}$ & $E^{ad}_6$ & $E^{ad}_7$ & $E^{ad}_8$ \\
\setlength{\unitlength}{2.5mm}
\begin{picture}(10,3)(-2.5,-1)
\multiput(0,0)(2,0){4}{$\circ$}
\put(0.55,.35){\line(1,0){1.55}}
\put(4.55,.35){\line(1,0){1.55}}
\multiput(2.55,.25)(0,.2){2}{\line(1,0){1.55}} \put(2.75,0){$>$}
\put(6,0){$\bullet$}
\put(0,0){$\ast$}
\end{picture} &
\setlength{\unitlength}{2.5mm}
\begin{picture}(10,3)(0,-1)
\multiput(0,0)(2,0){5}{$\circ$}
\multiput(0.55,.35)(2,0){4}{\line(1,0){1.55}}
\put(4.3,-1.4){\line(0,1){1.5}}
\put(4,-2){$\circ$}
\put(0,0){$\bullet$}
\put(8,0){$\bullet$}
\put(4,-2){$\ast$}
\end{picture} &
\setlength{\unitlength}{2.5mm}
\begin{picture}(11,3)(0.8,-1)
\multiput(0,0)(2,0){6}{$\circ$}
\multiput(0.5,.3)(2,0){5}{\line(1,0){1.6}} \put(4,-2){$\circ$}
\put(4.3,-1.4){\line(0,1){1.5}}
\put(0,0){$\ast$}
\put(8,0){$\bullet$}
\end{picture} &
\setlength{\unitlength}{2.5mm}
\begin{picture}(11,3)(0.8,-1)
\multiput(0,0)(2,0){7}{$\circ$}
\multiput(0.5,.3)(2,0){6}{\line(1,0){1.6}} \put(4,-2){$\circ$}
\put(4.3,-1.4){\line(0,1){1.5}}
\put(0,0){$\bullet$}
\put(12,0){$\ast$}
\end{picture} \\
\end{tabular}\end{center}

\section{Hyperplane sections}

Given an algebraic variety $X\subset\ppp V$, it is interesting to study  
the hyperplane sections $X\cap H$, for example the variation of
topology or Hodge structure
as one varies the hyperplane $H$. 
When $X=G/P$ is homogeneous one could hope to have
explicit descriptions of all hyperplane sections, at least when there
are a finite number of $G$-orbits on $V$. Donagi \cite{donagi} gives
such explicit descriptions for the Grassmanian  $G(3,6)$ and we
generalize his description to the F-planes of the third row 
of the chart, as well as recording the sections of the Severi
varieties.
 
\begin{prop} There are three types of hyperplane sections of a 
Severi variety:
\begin{enumerate}
\item homogeneous sections, this is the generic case;
\item sections with a unique singular point, which is an ordinary
  quadratic singularity;
\item sections whose singular locus is an $\BA\pp 1 (=Q^m)$.
\end{enumerate}
\end{prop}

\begin{prop}
The four types of hyperplane sections of $G_w(\BA^3, \BA^6)$ are:
\begin{enumerate} 
\item smooth generic sections;
\item sections with a unique singularity, an ordinary quadratic
singularity;
\item sections that are singular along a smooth quadric of dimension 
$m+1$;
\item sections whose singular locus is a cone over an $\BA\pp 2$.
\end{enumerate}
\end{prop}
 
These descriptions follow from our discussions above. More precisely, 
 if $G$ acts   on $V$ with a finite number of orbits
and with closed
orbit $X\subset\BP V$, the orbit structure of $G$
on $V^*$ is the same as that on $V$. 
A generic hyperplane section is always smooth, and if
the dual variety is nondegenerate, smooth points of the dual
variety give rise to sections with a unique singularity that
is an ordinary quadratic singularity. (In general, if $X^*$ has defect
$\dstar$, and $H\in X^*_{smooth}$,
 then $(X\cap  H)_{sing}$ is a $\pp{\dstar}$.)

In the Severi case,  the third type of section follows from Freudenthal's
perspective, an F-line has contact with an $\BA\pp 1$.  
In the $G_w(\BA^3,\BA^6)$
case, the third and fourth types of sections again follow  from Freudenthal
where we identify $[y]\in \ppp V$ with $[y^*]= [\Omega(y,\cdot )]\in
\ppp V^*$. We have $H\supset \tilde T_zG_w(\BA^3,\BA^6)$ if and only
if $H^*\in \tilde T_zG_w(\BA^3,\BA^6)$
and $H$ and $H^*$ are in isomorphic orbit closures.   The vertex of the
cone in
the last case
 is  of course $H^*$.

\end{document}